\documentclass[11pt, a4paper,twoside]{article}
\usepackage{amssymb,amsfonts,amsmath,amsthm,euscript}
\usepackage[english]{babel}
\usepackage[latin1]{inputenc}
\usepackage{dsfont}
\usepackage{geometry}
\usepackage{graphicx}
\usepackage{mathrsfs}
\usepackage{multirow}
\usepackage{multicol}
\usepackage[bf,footnotesize]{caption}
\usepackage{sectsty}
\usepackage{mathtools}
\usepackage{setspace}
\usepackage{placeins}
\usepackage{subfigure}
\sectionfont{\large}
\allowdisplaybreaks 
\DeclareMathAlphabet{\mz}{OT1}{pzc}{m}{it}
\DeclareMathOperator*{\diag}{\mathrm{diag}}
\DeclareMathOperator*{\tr}{\mathrm{tr}}

\begin{document}

\theoremstyle{plain}
\newtheorem{thm}{Theorem}[section]
\newtheorem{cor}{Corollary}[section]
\newtheorem{prop}{Proposition}[section]
\newtheorem{lema}{Lemma}[section]

\theoremstyle{definition}
\newtheorem{rmk}{Remark}[section]

\renewcommand{\theequation}{\thesection.\arabic{equation}}
\addtolength{\parskip}{.2cm}
\newcommand{\comp}[1]{\overline{#1}^{\,\,\prime}}
\newcommand{\re}{\mathrm{Re}}
\newcommand{\pim}{\mathrm{Im}}
\newcommand{\R}{\mathds{R}}
\newcommand{\N}{\mathds{N}}
\newcommand{\E}{\mathds{E}}
\newcommand{\Z}{\mathds{Z}}
\newcommand{\C}{\mathds{C}}
\renewcommand{\P}{\mathds{P}}
\renewcommand{\l}{\lambda}
\renewcommand{\L}{\Lambda}
\renewcommand{\d}{\widehat{\bs d}}
\renewcommand{\do}{\bs d_0}
\newcommand{\G}[1]{\widehat G(#1)}
\newcommand{\od}{\overline{\bs {\mathrm{d}}}}
\newcommand{\sz}{\frac{1}{m}\sum_{j=[m\kappa]}^m}
\newcommand{\qz}{\mathcal{Q}_{\kappa}(\bs d)}
\newcommand{\eps}{\varepsilon}
\newcommand{\te}{\theta}
\newcommand{\mcd}{\widehat{\mathcal{D}}}
\newcommand{\mdk}{\widehat{\mathcal{D}}_{\kappa}}
\newcommand{\mcp}{\mathcal{P}}
\newcommand{\im}{\mathrm{i}}
\newcommand{\sk}{\sum_{k=1}^q}
\newcommand{\bs}[1]{\boldsymbol{#1}}
\renewcommand{\proof}{ \noindent \textbf{\emph{Proof: }}}
\newcommand{\fim}{\hfill{\footnotesize$\blacksquare$}\normalsize\\}
\newcommand{\id}{\mathrm{I}}
\newcommand{\lfl}{\L_j(\do)^{-1}f_n(\l_j)\comp{\L_j(\do)^{-1}}}
\newcommand{\ir}{\mathrm{I}_{(r)}}
\newcommand{\is}{\mathrm{I}_{(s)}}
\numberwithin{equation}{section}
\numberwithin{table}{section}
\numberwithin{figure}{section}
\pagestyle{myheadings} 
\markboth{A Generalization of a Gaussian Semiparametric Estimator}{G. Pumi and S.R.C. Lopes} 
\long\def\sfootnote[#1]#2{\begingroup%
\def\thefootnote{\fnsymbol{footnote}}\footnote[#1]{#2}\endgroup}

\def\bfootnote{\xdef\@thefnmark{}\@footnotetext}

\thispagestyle{empty}
\vskip2cm
{\centering
\Large{\bf  A Generalization of a Gaussian Semiparametric Estimator on Multivariate Long-Range Dependent Processes}\vspace{.6cm}\\
\large{ {\bf Guilherme Pumi$\!\phantom{i}^{\mathrm{a,}}$\sfootnote[1]{Corresponding author.}\let\thefootnote\relax\footnote{\hskip-.3cm$\phantom{s}^\mathrm{a}$Mathematics Institute - Federal University of Rio Grande do Sul - 9500,  Bento Gon\c calves Avenue - 91509-900, Porto Alegre - RS - Brazil.}\let\thefootnote\relax\footnote{E-mail addresses: guipumi@gmail.com (G. Pumi), silvia.lopes@ufrgs.br (S.R.C. Lopes).} and S\'ilvia R.C. Lopes$\!\!\phantom{s}^\mathrm{a}$ } \\
\let\thefootnote\relax\footnote{This version: \today.}}}
\vspace{-.2cm}
\begin{center} \large{
Mathematics Institute\\
Federal University of Rio Grande do Sul}
\end{center}

\vskip.6cm

\begin{abstract}
In this paper we propose and study a general class of Gaussian Semiparametric Estimators (GSE) of the fractional differencing parameter in the context of long-range dependent multivariate time series. We establish large sample properties of the estimator without assuming Gaussianity. The class of models considered here satisfies simple conditions on the spectral density function, restricted to a small neighborhood of the zero frequency and  includes important class of VARFIMA processes.  We also present a simulation study to assess the finite sample properties of the proposed estimator based on a smoothed version of the GSE which supports its competitiveness.\vspace{.2cm}\vspace{.2cm}
\end{abstract}

\section{Introduction}
Let $\bs d=(d_1,\cdots,d_q)^\prime\in(-1/2,1/2)^q$ and let  $\cal B$ be the shift operator. Consider the $q$-dimensional weakly stationary process $\{\bs{X}_{\!t}\}_{t=0}^\infty$ obtained as a stationary solution of the difference equations
\small\begin{equation}\label{specif}
\diag_{k\in\{1,\cdots,q\}}\!\!\!\big\{(1-\mathcal{ B})^{d_k}\!\big\}\big(\bs{X}_{\!t}-\E(\bs{X}_{\!t})\big)=\bs{Y}_{\!\!t},
\end{equation}\normalsize
where $\{\bs{Y}_{\!\!t}\}_{t=0}^\infty$  is a $q$-dimensional weakly stationary process whose spectral density function $f_{\bs Y}$ is bounded and bounded away from zero.  Each coordinate process in \eqref{specif} exhibits long-range dependence whenever the respective parameter $d_i>0$, in the sense that the spectral density function satisfies $f(\l)\sim K\l^{-2d_i}$, as $\l\rightarrow0^+$, for  some constant $K>0$ and $i\in \{1,\cdots, q\}$.

Processes of the form \eqref{specif} constitute the so-called fractionally integrated processes. As a particular case, consider the situation where the $i$-th coordinate process \small $\big\{Y_t^{(i)}\big\}_{t=0}^\infty$ \normalsize follows an ARMA model. In this case, the associated coordinate process \small $\big\{X_t^{(i)}\big\}_{t=0}^\infty$ \normalsize will be a classical ARFIMA process with the same AR and MA orders and differencing parameter $d_i$. If the process $\{\bs{Y}_{\!\!t}\}_{t=0}^\infty$ is a vectorial  ARMA process, then the resulting multivariate process will be the so-called VARFIMA process with differencing parameter $\bs d=(d_1,\cdots,d_q)^\prime$.  VARFIMA and, more generally, fractionally integrated processes, are widely used to model multivariate processes with long-range dependence. See, for instance, the recent work of Chiriac and Voev (2011) on modeling and forecasting high frequency data by using VARFIMA and fractionally integrated processes.

The parameter $\bs d$ in \eqref{specif} determines the spectral density function behavior at the zero frequency as well as the long run  autocovariance/autocorrelation structure. Hence, estimation becomes an important matter whenever the long run structure of the process is of interest.

Estimation of the parameter $\bs d$ in the multivariate case has seen a growing interest in the last years. A maximum likelihood approach was first considered in Sowell (1989), but the computational cost of the author's method at the time was very high A few years later, Luce\~no (1996) presented a computationally cheaper alternative for the maximum likelihood approach based on rewriting and approximating the quadratic form of the Gaussian likelihood function. In a recent work, Tsay (2010) proposed an even faster approach to calculate the exact conditional likelihood based on the multivariate Durbin-Levinson algorithm. Although the maximum likelihood approach usually  provides good results, it is still a computationally expensive method.

The works of Fox and Taqqu (1986), Giraitis and Surgailis (1990), among others, provided a rigorous asymptotic theory for (univariate) Gaussian parametric estimates which includes, for instance, $n^{1/2}$-consistency and asymptotic normality. One drawback is the crucial role played by the Gaussianity assumption in the theory, which also requires strong distributional and regularity conditions and is non-robust with respect to the parametric specification of the model, leading to inconsistent estimates under misspecification.

In the univariate case, Gaussian Semiparametric Estimation (GSE) was first introduced in K\"unsch (1987) and later rigorously developed by Robinson (1995b). It provides a more robust alternative compared to the parametric one, requiring less distributional assumptions and being more efficient. In the multivariate case, Robinson (1995a) was the first to study and develop a rigorous treatment of a semiparametric estimator. A two-step multivariate GSE has been studied in the work of Lobato (1999), which showed its asymptotic normality under mild conditions, but without relying on Gaussianity. A few years later, Shimotsu (2007) introduced a refinement of Lobato's two-step GSE, which is consistent and asymptotically normal under very mild conditions (Gaussianity is, again, nowhere assumed), but with smaller asymptotic variance than Lobato's estimator. The technique applied in Shimotsu (2007) was a multivariate extension of that in Robinson (1995b), powerful enough to show not only the consistency of the proposed estimator, but also the consistency of Lobato's two-step GSE. Recently, Nielsen (2011) extended the work of Shimotsu (2007) to include the non-stationary case by using the so-called extended periodogram while Pumi and Lopes (2013) extends the work of Lobato (1999) by considering general estimators of the spectral density function in Lobato's objective function.

The estimator introduced in Shimotsu (2007) is based on the specification of the spectral density function in a neighborhood of the zero frequency. Estimation of the differencing parameter $\bs d$ is obtained through minimization of an objective function, which is derived from the expression of the Gaussian log-likelihood function near the zero frequency. To obtain the objective function, the spectral density is estimated by the periodogram of the process. Although asymptotically unbiased, it is well known that the periodogram is not a consistent estimator of the spectral density, presents wild fluctuations near the zero frequency and, understood as a sequence of random variables, it does not converge to a random variable at all (cf. Grenander, 1951). Some authors actually consider the periodogram ``\emph{an extremely poor (if not useless) estimate of the spectral density function}'' (Priestley, 1981, page 420). The estimators introduced in Lobato (1999) and Shimotsu (2007) are known to present a good finite sample performance, but given the wild behavior of the periodogram, a natural question is: can we do better with a better behaved spectral density estimator? In this work, our primary goal is to provide an answer to this question.

Our contribution to the theory of GSE is two-folded. First, being consistency a highly desirable property of an estimator, we study the consequences of substituting the periodogram in Shimotsu (2007)'s objective function by an arbitrary consistent estimator of the spectral density function. We prove the consistency of the proposed estimator under the same assumptions as in Shimotsu (2007) and no assumption on the spectral density estimator other than consistency. Second, considering Shimotsu (2007)'s objective function with the periodogram substituted by an arbitrary spectral density estimator,  we derive necessary conditions under which GSE is consistent and satisfy a multivariate  CLT. Gaussianity is nowhere assumed. In order to assess the finite sample properties of the estimators studied here and its competitiveness, we present a simulation study based on simulated VARFIMA process. We apply the so-called smoothed periodogram and the tapered periodogram as estimators of the spectral density function.

The paper is organized as follows. In the next section, we present some preliminaries concepts and results necessary for this work and introduce a general class of estimators based on appropriate modifications of the Shimotsu's objective function. Section 3 is devoted to derive the consistency of the proposed estimator while in Section 4 we derive conditions for the proposed estimator to satisfy a multivariate CLT. In Section 5 we present some Monte Carlo simulation results to assess the finite sample performance of the proposed estimator. Conclusions and final remarks are reserved to Section 6. The usually long proofs of our results are postponed to the Appendix A.

\section{Preliminaries}

Let $\{\bs{X}_{\!t}\}_{t=0}^\infty$ be a $q$-dimensional process specified by \eqref{specif} and assume that the spectral density matrix of $\bs{Y}_{\!\!t}$ satisfies $f_{\bs Y}\sim G$ for a real, symmetric, finite and positive definite matrix $G$. Let $f$ be the spectral density matrix function of $\bs{X}_{\!t}$, so that
\small\[\E\big[\big(\bs{X}_{\!t}-\E(\bs{X}_{\!t})\big)\big(\bs{X}_{\!t+h}-\E(\bs{X}_{\!t})\big)'\big]=\int_{-\pi}^\pi \mathrm{e}^{\im h \l}f(\l)d\l,\]\normalsize
for $h\in\N^\ast:=\N\!\setminus\!\{0\}$. Following the reasoning in Shimotsu (2007), the spectral density matrix of $\bs{X}_{\!t}$ at the Fourier frequencies $\l_j=2\pi j/n$, with $j=1,\cdots,m$ and $m=o(n)$, can be written as
\small\begin{equation}\label{fapprox}
f(\l_j)\sim \L_j(\bs d)G\comp{\L_j(\bs d)},\quad \mbox{for} \quad \L_j(\bs d)=\diag_{k\in\{1,\cdots,q\}}\{\L_j^{(k)}(\bs d)\}\quad\mbox{and} \quad \L_j^{(k)}(\bs d)=\l_j^{-d_k}\mathrm{e}^{\im(\pi-\l_j)d_k/2},
\end{equation}\normalsize
where, for a complex matrix $A$, $\comp{A}$ denotes the conjugate transpose of $A$.  Let
\small\[I_n(\l)\vcentcolon=w_n(\l)\comp{w_n(\l)},\quad\mbox{ where }\quad w_n(\l)\vcentcolon=\frac{1}{\sqrt{2\pi n}}\sum_{t=1}^n\bs X_{\!t} \mathrm{e}^{\im t\l},\]\normalsize
be the periodogram and the discrete Fourier transform of $\bs X_{\! t}$ at $\l$, respectively.
From the local form of the spectral density at zero frequency, given in \eqref{fapprox}, replaced in the frequency domain Gaussian log-likelihood localized at the origin, Shimotsu (2007) proposed a semiparametric estimator for the fractional differencing parameter $\bs d$ based on the objective function
\small\begin{equation}\label{shimR}
R(\bs d)\vcentcolon=\log\big(\det \{\tilde G(\bs d)\}\big)-2\sum_{k=1}^qd_k\frac{1}{m}\sum_{j=1}^m\log(\l_j),
\end{equation}\normalsize
where
\small\begin{equation}\label{shimG}
\tilde G(\bs d)\vcentcolon=\frac{1}{m}\sum_{j=1}^m\re\big[\L_j(\bs d)^{-1}I_n(\l_j)\comp{\L_j(\bs d)^{-1}}\big].
\end{equation}\normalsize
with $\L_j(\bs d)$ defined in \eqref{fapprox}. The estimator of $\bs d$ is then given by
\small\begin{equation}\label{shimest}
\widetilde{\bs d}=\underset{\bs d\in\Theta}{\arg\min}\{R(\bs d)\},
\end{equation}\normalsize
where the space of admissible estimates is of the form $\Theta=[-1/2+\epsilon_1,1/2-\epsilon_2]$, for arbitrarily small $\epsilon_i>0$, $i=1,2$, henceforth fixed except stated otherwise. Shimotsu (2007) shows that the estimator based on the objective function \eqref{shimR} and \eqref{shimG} is consistent under mild conditions. Given the wild behavior of the periodogram as an estimator of the spectral density function, specially near the origin, in this work we consider substituting the periodogram by some other spectral density estimator, say $f_n$. Our interest lies on estimators based on objective functions of the form
\small\begin{equation}\label{sshR}
S(\bs d)\vcentcolon=\log\big(\det \{\widehat G(\bs d)\}\big)-2\sum_{k=1}^qd_k\frac{1}{m}\sum_{j=1}^m\log(\l_j),
\end{equation}\normalsize
 with
\small\begin{equation}\label{sshG}
\widehat G(\bs d)\vcentcolon=\frac{1}{m}\sum_{j=1}^m\re\big[\L_j(\bs d)^{-1}f_n(\l_j)\comp{\L_j(\bs d)^{-1}}\big].
\end{equation}\normalsize
Notice that \eqref{sshR} is just \eqref{shimR} with the periodogram $I_n$ in \eqref{shimG} replaced by $f_n$. The estimator of $\bs d$ is then defined analogously as
\small\begin{equation}\label{estimator}
\widehat{\bs d}=\underset{\bs d\in\Theta}{\arg\min}\{S(\bs d)\}.
\end{equation}\normalsize
In the sections to come, we shall study the asymptotic behavior of estimator \eqref{estimator}. The study is focused on two different classes of spectral density estimator. First we shall consider the class of consistent estimators of the spectral density function and we shall show that, under no further hypothesis on the $f_n$, estimator \eqref{estimator} is consistent. Second, we consider a class of spectral density functions satisfying a moment condition and we shall derive conditions for the consistency and asymptotic normality of the estimator \eqref{estimator}.

Before proceeding with the results, we shall establish some notation. Let $\{\bs{X}_{\!t}\}_{t=0}^\infty$ be a $q$-dimensional process specified by \eqref{specif} and let $\{\bs \eps_t\}_{t\in\Z}$ such that $\bs{X}_{\!t}-\E(\bs{X}_{\!t})=\sum_{k=0}^\infty A_k\bs\eps_{t-k}$.
We define a function $A$ by setting
\small\begin{equation}\label{ast}
A(\l)\vcentcolon=\sum_{k=0}^\infty A_k\mathrm{e}^{\im k\l}.
\end{equation}\normalsize
The periodogram function associated to $\{\bs \eps_t\}_{t\in\Z}$ is denoted by $I_{\bs\eps}$, that is,
\small\begin{equation}\label{per_eps}
I_{\bs \eps}(\l)\vcentcolon=w_{\bs\eps}(\l)\comp{w_{\bs\eps}(\l)},\qquad\text{where}\qquad w_{\bs\eps}(\l)\vcentcolon=\frac{1}{\sqrt{2\pi n}}\sum_{t=1}^n\bs\eps_t\mathrm{e}^{\im t\lambda}.
\end{equation}\normalsize
For a matrix $M$, we shall denote the $r$-th row and the $s$-th column of $M$ by $(M)_{r\bs\cdot}$ and $(M)_{\bs\cdot s}$, respectively.

\section{Consistency of the estimator}

Let $\{\bs{X}_{\!t}\}_{t=0}^\infty$ be a $q$-dimensional process specified by \eqref{specif} and $f$ be its spectral density matrix. Suppose that the spectral density matrix of the weakly stationary process $\{\bs{Y}_{\!\!t}\}_{t=0}^\infty$ in \eqref{specif} satisfies $f_{\bs Y}(\l)\sim G_0$ for a real, symmetric and positive definite matrix $G_0=(G_0^{rs})_{r,s=1}^q$. Let $\bs d_0=(d_1^0,\cdots,d_q^0)^\prime$ be the true fractional differencing vector parameter and assume that the following assumptions are satisfied:
\begin{itemize}
\item[\textbf{A}1.] As $\l\rightarrow0^+$,
\small\[f_{rs}(\lambda)=\mathrm{e}^{\im\pi(d_r^0-d_s^0)/2}G_0^{rs}\l^{-d_r^0-d_s^0}+o(\l^{-d_r^0-d_s^0}),\quad \mbox{ for all }\,r,s=1,\cdots,q.\]\normalsize
\item[\textbf{A}2.] Denoting the sup-norm by $\|\cdot\|_{\infty}$, assume that
\small\begin{equation}\label{condema}
\bs{X}_{\!t}-\E(\bs{X}_{\!t})=\sum_{k=0}^\infty A_k\bs\eps_{t-k},\quad\sum_{k=0}^\infty\big\|A_k\big\|_{\infty}^2<\infty,
\end{equation}\normalsize
where $\{\bs \eps_t\}_{t\in\Z}$ is a process such that
\small\[\E(\bs\eps_t|\mathscr{F}_{t-1})=0\quad \mbox{ and }\quad\E(\bs\eps_t\bs\eps_t^\prime|\mathscr F_{t-1})=\id_q,\quad\mbox{a.s.}\]\normalsize
for all $t\in\Z$, where $\id_q$ is the $q\times q$ identity matrix and $\mathscr{F}_t$ denotes the $\sigma$-field generated by $\{\bs\eps_s, s\leq t\}$. Also assume that there exist a scalar random variable $\xi$ and a constant $K>0$ such that $\E(\xi^2)<\infty$ and $\P\big(\|\bs\eps_t\|_\infty^2>\eta\big)\leq K\P(\xi^2>\eta)$, for all $\eta>0$.
\item[\textbf{A}3.] In a neighborhood $(0,\delta)$ of the origin, $A$ given by \eqref{ast} is differentiable and, as $\l\rightarrow0^+$,
\small\[\frac{\partial}{\partial\l}\big(\comp{A(\l)}\big)_{r\bs\cdot }=O\big(\l^{-1}\big\|\big(\comp{A(\l)}\big)_{r\bs\cdot }\big\|_{\infty}\big).\]\normalsize
\item[\textbf{A}4.] As $n\rightarrow\infty$,
\small\[\frac{1}{m}+\frac{m}{n}\longrightarrow 0.\]\normalsize
\end{itemize}

\begin{rmk}
Assumptions \textbf{A}1-\textbf{A}4 are the same as in Shimotsu (2007) and are multivariate extensions of the assumptions made in Robinson (1995b) and analogous to the ones used in Robinson (1995a) and Lobato (1999). Assumption \textbf{A}1 describes the true spectral density matrix behavior at the origin. Notice that, since $\lim_{\l\rightarrow 0^+}\mathrm{e}^{\im \lambda}-1=0$, replacing $\mathrm{e}^{\im\pi(d_r^0-d_s^0)/2}$ by $\mathrm{e}^{\im(\pi-\lambda)(d_r^0-d_s^0)/2}$ makes no difference. Assumption \textbf{A}2 regards the causal representation of $\bs X_t$, and more specifically, the behavior of the innovation process which is assumed to be a not necessarily uncorrelated square integrable martingale difference  uniformly dominated (in probability) by a scalar random variable with finite second moment. Assumption \textbf{A}3 is a regularity condition (also imposed in Fox and Taqqu, 1986 and Giraitis and Surgailis, 1990, among others, in the parametric case) and will be useful in proving Lemmas \ref{lema1} and \ref{lema2} below. Assumption \textbf{A}4 is minimal but necessary since $m$ must go to infinity for consistency, but slower than $n$ in view of Assumption \textbf{A}1.
\end{rmk}

Observe that assumptions \textbf A1-\textbf{A}4 are only concerned to the behavior of the spectral density matrix on a neighborhood of the origin and, apart from integrability (implied by the process's weakly stationarity property), no assumption whatsoever is made on the spectral density matrix behavior outside this neighborhood. For $\beta\in(0,1)$, let $f_n$ be an $n^{\beta}$-consistent (that is, $n^{\beta}(f_n-f)\overset{\P}{\longrightarrow}0$, as $n$ goes to infinity) estimator of the spectral density for all $\do\in B$, where $B\subset\R^q$ is a closed set. If $d_k^{0}\in(0,0.5)$, the respective component of the spectral density matrix of $\{\bs X_t\}_{t=0}^\infty$ is unbounded at the origin. Hence, there is no hope in obtaining a consistent estimator of the spectral density function when $\do\in(0,0.5)^q$. For $q\in\N^\ast$, let

\small\begin{equation}\label{set}
\Omega_\beta\vcentcolon=\bigg[-\frac{\beta}{2},0\bigg)^q\bigcap\bigg(-\frac{1}{2},0\bigg)^q\bigcap B\subseteq\bigg(-\frac12,0\bigg)^q.
\end{equation}\normalsize
Lemma \ref{lema1} establishes the consistency of $\widehat G(\do)$ given in \eqref{sshG} under the assumption of $n^\beta$-consistency of $f_n$ in $B$. Due to their lengths, the proofs of all results in the paper are postponed to the Appendix A.

        \begin{lema}\label{lema1}
        Let $\{\bs{X}_{\!t}\}_{t=0}^\infty$ be a $q$-dimensional process specified by \eqref{specif} and $f$ be its spectral density matrix. Let $f_n$ be a $n^\beta$-consistent estimator for $f$, for all $\do\in B$. If $\bs d_0\in\Omega_\beta$, then
        \small\[\widehat G(\do)=G_0+o_\P(1).\]\normalsize
        \end{lema}

Theorem \ref{consistency} establishes the consistency of $\d$, given in \eqref{estimator} with $\Theta$ substituted by $\Omega_\beta$, under assumptions  \textbf{A}1-\textbf{A}4 and assuming $n^\beta$-consistency of the spectral density function estimator.

        \begin{thm}\label{consistency}
        Let $\{\bs{X}_{\!t}\}_{t=0}^\infty$ be a $q$-dimensional process specified by \eqref{specif} and $f$ be its spectral density matrix. Let $f_n$ be a $n^\beta$-consistent estimator of $f$, for all $\do\in B$ and $\beta\in(0,1)$, and let $\widehat{\bs d}$ be as in \eqref{estimator} with $\Theta$ substituted by $\Omega_\beta$.  Assume that assumptions \textbf{A}1-\textbf{A}4 hold and let $\bs d_0\in\Omega_\beta$. Then, $\widehat{\bs d}\overset{\P}{\longrightarrow}\bs d_0$, as $n\rightarrow \infty$.
        \end{thm}

Assuming the consistency of the spectral density estimator $f_n$ in Theorem \ref{consistency} excluded the case $\bs d_0\in(0,0.5)^q$, so that, under this assumption,  the process $\{\bs{X}_{\!t}\}_{t=0}^\infty$ can have no long-range dependent component. To overcome this limitation, we now consider the class $\mathscr D$ of estimators $f_n:=(f_n^{rs})_{r,s=1}^q$ satisfying, for all $r,s\in\{1,\cdots,q\}$,
\small\begin{equation}\label{condthm}
\E\Big(\l_j^{d_r^0+d_s^0}\Big|f_n^{rs}(\l_j)-\big(A(\l_j)\big)_{r\bs\cdot}I_{\bs\eps}(\l_j)\big(\comp{A(\l_j)}\big)_{\bs\cdot s}\Big|\Big)\:= o(1), \quad \mbox{ as }\:n \rightarrow \infty,
\end{equation}\normalsize
where $A$ and $I_{\bs\eps}$ are given by \eqref{ast} and \eqref{per_eps}, respectively, and $\bs d_0\in\Theta\subset[-0.5,0.5]^q$. Condition \eqref{condthm} is satisfied by the ordinary periodogram and the tapered periodogram and, thus, $\mathscr D$ is non-empty. The next lemma will be useful in proving Theorem \ref{nonbeta}.

        \begin{lema}\label{lema2}
        Let $\{\bs{X}_{\!t}\}_{t=0}^\infty$ be a $q$-dimensional process specified by \eqref{specif} and $f$ be its spectral density matrix. Let $f_n\in\mathscr D$ and assume that assumptions \textbf{A}1-\textbf{A}4 hold. Then, for $1\leq u<v\leq m$,
        \small\[\max_{r,s\in\{1,\cdots,q\}}\bigg\{\sum_{j=u}^v\mathrm{e}^{\im(\l_j-\pi)(d_r^0-d_s^0)/2}\l_j^{d^0_r+d^0_s}f_n^{rs}(\l_j)-G_0^{rs}\bigg\}= \mathscr{A}_{uv} +\mathscr{B}_{uv},\]\normalsize
        where $\mathscr A_{uv}$ and $\mathscr B_{uv}$ satisfy
        \small\[\E\big(|\mathscr A_{uv}|\big)=o(v-u+1)\qquad\mbox{ and }\qquad \max_{1\leq u<v\leq m}\big\{\big|v^{-1}\mathscr{B}_{uv}\big|\big\}=o_\P(1).\]\normalsize
        \end{lema}

In Theorem \ref{nonbeta} we derive a necessary condition for the consistency of $\d$ given in \eqref{estimator}, when the consistency condition on $f_n$ is relaxed and we assume $f_n\in\mathscr D$ instead.

        \begin{thm}\label{nonbeta}
        Let $\{\bs{X}_{\!t}\}_{t=0}^\infty$ be a $q$-dimensional process specified by \eqref{specif} and $f$ be its spectral density matrix. Let $f_n\in\mathscr D$ be an estimator of $f$, and consider the estimator $\widehat{\bs d}$, based on $f_n$, given in \eqref{estimator}.  Assume that assumptions \textbf{A}1-\textbf{A}4 hold. Then, $\widehat{\bs d}\overset{\P}{-\!\!\!\longrightarrow}\bs d_0$, as $n\rightarrow \infty$.
        \end{thm}

\section{Asymptotic normality of the estimator}
In this section we present a sufficient condition for the asymptotic normality of the GSE given by \eqref{estimator}, under similar assumptions as in Shimotsu (2007), with $f_n$ an estimator of the spectral density function satisfying a single regularity condition. The asymptotic distribution of the estimator \eqref{estimator} will be the same as \eqref{shimest}, established by Shimotsu (2007). 

Again, let $\{\bs{X}_{\!t}\}_{t=0}^\infty$ be a $q$-dimensional process specified by \eqref{specif} and $f$ be its spectral density matrix. Suppose that the spectral density matrix of the weakly stationary process $\{\bs{Y}_{\!\!t}\}_{t=0}^\infty$ in \eqref{specif} satisfies $f_{\bs Y}(\l)\sim G_0$ for a real, symmetric and positive definite matrix $G_0=(G_0^{rs})_{r,s=1}^q$. Let $\bs d_0=(d_1^0,\cdots,d_q^0)^\prime$ be the true fractional differencing vector parameter. Assume that the following assumptions are satisfied
\begin{itemize}
\item[\textbf B1.] For $\alpha\in(0,2\,]$ and $r,s\in\{1,\cdots,q\}$,
\small\[f_{rs}(\l)=\mathrm{e}^{\im(\pi-\l)(d_r^0-d_s^0)/2}\l^{-d_r^0-d_s^0}G_0^{rs}+O(\l^{-d_r^0-d_s^0+\alpha}),\,\,\text{ as }\l\rightarrow0^+.\]\normalsize
\item[\textbf B2.] Assumption \textbf A2 holds and the process $\{\bs\eps_t\}_{t\in\Z}$ has finite fourth moment.
\item[\textbf B3.] Assumption \textbf A3 holds.
\item[\textbf B4.] For any $\delta>0$,
\small\[\frac{1}{m}+\frac{m^{1+2\alpha}\log(m)^2}{n^{2\alpha}}+\frac{\log(n)}{m^\delta}\longrightarrow 0, \,\,\text{ as } n\rightarrow\infty.\]\normalsize
\item[\textbf B5.] There exists a finite real matrix $M$ such that
\small\[\L_j(\do)^{-1}A(\l_j)=M+o(1),\,\,\text{ as }\l_j\rightarrow0.\]\normalsize
\end{itemize}
\vskip.2cm
\begin{rmk}\label{rmk}
Assumption \textbf B1 is a smoothness condition often imposed in spectral analysis. Compared to assumption 1 in Robinson (1995a), assumption  \textbf B1 is slightly more restrictive. It is satisfied by certain VARFIMA processes. Assumption \textbf B2 imposes that the process $\{\bs X_t\}_{t\in\N^\ast}$ is linear with finite fourth moment. This restriction in the innovation process is necessary since at a certain point of the proof of the asymptotic normality, we shall need a CLT result regarding a certain martingale difference derived from a quadratic form involving  $\{\bs\eps_t\}_{t\in\Z}$, which must have finite variance. Assumption \textbf B4 is the same as assumption $4'$ in Shimotsu (2007) and is slightly stronger than the ones imposed in Robinson (1995b) and Lobato (1999) (see Shimotsu, 2007 p.283 for a discussion). It implies that $(m/n)^b=o\big(m^{-\frac{b}{2\alpha}}\,\log(m)^{-\frac{b}{\alpha}}\big)$, for $b\neq 0$.  Assumption \textbf B5 is the same as assumption 5' in Shimotsu (2007) and is a mild regularity condition in the degree of approximation of $A(\l_j)$ by $\L_j(\do)$.  It is satisfied by general VARFIMA processes.
\end{rmk}

The next lemma will be useful in proving Theorem \ref{norm}. The proofs of the results in this section are presented in Appendix A.

        \begin{lema}\label{lemma1b2}
         Let $\{\bs{X}_{\!t}\}_{t=0}^\infty$ be a $q$-dimensional process specified by \eqref{specif} and $f$ be its spectral density matrix. Let $f_n$ be an estimator of $f$, and consider the estimator $\widehat{\bs d}$, based on $f_n$, given in \eqref{estimator}.  Assume that assumptions \textbf{B}1-\textbf{B}5 hold and that $f_n$ satisfies
         \small\begin{equation}\label{cond_an}
        \max_{1\leq v\leq m}\bigg\{\sum_{j=1}^v\big[f_n^{rs}(\l_j)-\big(A(\l_j)\big)_{r\bs\cdot}I_{\bs\eps}(\l_j)\big(\comp{A(\l_j)}\big)_{\bs\cdot s}\big]\bigg\}=o_\P\bigg(\frac m {n^{1+|d_r^0+d_s^0|}}\bigg),
         \end{equation}\normalsize
         for all $r,s\in\{1,\cdots,q\}$ and $\do\in\Theta$, where $A$ and $I_{\bs\eps}$ are defined in \eqref{ast} and \eqref{per_eps}, respectively. Then,
         \begin{itemize}
        \item[\rm {(a)}] uniformly in $1\leq v \leq m$,
        \small\begin{equation}\label{lema1b2a}
        \hspace{-.7cm}\max_{r,s\in\{1,\cdots,q\}}\!\!\bigg\{\sum_{j=1}^v\mathrm{e}^{\im(\l_j-\pi)(d_r^0-d_s^0)/2}\l_j^{d_r^0+d_s^0}\Big[f_n^{rs}(\l_j)- \big(A(\l_j)\big)_{r\bs\cdot}I_{\bs\eps}(\l_j)\big(\comp{A(\l_j)}\big)_{\bs\cdot s}\Big]\bigg\}=o_\P\bigg(\frac{m^{1/2}}{\log(m)}\bigg);
        \end{equation}\normalsize
        \item[\rm {(b)}] uniformly in $1\leq v \leq m$,
        \small\begin{equation}\label{lema1b2b}
        \max_{r,s\in\{1,\cdots,q\}}\!\!\bigg\{\sum_{j=1}^v\mathrm{e}^{\im(\l_j-\pi)(d_r^0-d_s^0)/2}\l_j^{d_r^0+d_s^0}f_n^{rs}(\l_j)-G_0^{rs}\bigg\} = O_\P\bigg(\frac{m^{\alpha+1}}{n^\alpha}+m^{1/2}\log(m)\bigg).
        \end{equation}\normalsize
         \end{itemize}
        \end{lema}

The next theorem presents a necessary condition for the asymptotic normality of the GSE given in \eqref{estimator}. We notice that the variance-covariance matrix of the limiting distribution is the same as the estimator in \eqref{shimest}, as derived in Shimotsu (2007).

        \begin{thm}\label{norm}
        Let $\{\bs{X}_{\!t}\}_{t=0}^\infty$ be a $q$-dimensional process specified by \eqref{specif} and $f$ be its spectral density matrix. Let $f_n$ be an estimator of $f$, and consider the estimator $\widehat{\bs d}$, based on $f_n$, given in \eqref{estimator}.  Assume that assumptions \textbf{B}1-\textbf{B}5 hold. Suppose that $f_n$ satisfies \eqref{cond_an},
        for all $r,s\in\{1,\cdots,q\}$ and $\bs d_0\in\Theta$. If $\d\overset{\P}{\longrightarrow}\do$, for $\do\in\Theta$, then
        \small\[m^{1/2}(\d-\do)\overset{d}{-\!\!\!\longrightarrow} N(\bs 0, \Sigma^{-1}),\]\normalsize
        as $n$ tends to infinity, where
        \small\[\Sigma:=2\bigg[G_0\odot G_0^{-1}+\id_q+\frac{\pi^2}{4}\big(G_0\odot G_0^{-1}-\id_q\big)\bigg],\]\normalsize
        with $\id_q$ the $q\times q$ identity matrix and $\odot$ denotes the Hadamard product.
        \end{thm}

\section{Monte Carlo Simulation Study }

In this section we perform a Monte Carlo simulation study to assess the finite sample performance of the estimator proposed in \eqref{estimator}.
We apply, as spectral density estimators, the so-called smoothed and tapered periodogram and compare them to the estimator \eqref{shimest}. We start recalling some facts about the smoothed and tapered periodogram.

\subsection{The Smoothed Periodogram}

Let $\{\bs{X}_{\!t}\}_{t=0}^\infty$ be a $q$-dimensional process specified by \eqref{specif}. Under some mild conditions, a class of consistent estimators of the spectral density of $\bs{X}_{\!t}$ is the so-called class of smoothed periodogram. For an array of functions $W_ n(k):=\big(W_n^{ij}(k)\big)_{i,j=1}^q$ (called weights) and $\{\ell(k)\}_{k=0}^\infty$ an increasing sequence of positive integers, the smoothed periodogram of $\{\bs{X}_{\!t}\}_{t=0}^\infty$ at the Fourier frequency $\l_j$ is defined as
\small\begin{equation}\label{smoothp}
\hat f_n(\l_j)\vcentcolon=\sum_{|k|\leq \ell(n)}W_n(k)\odot w_n(\l_{j+k})\comp{w_n(\l_{j+k})},
\end{equation}\normalsize
where $\odot$ denotes the Hadamard product.

The smoothed periodogram \eqref{smoothp} is a multivariate extension of the univariate smoothed periodogram. Notice that the use of the Hadamard product in \eqref{smoothp} allows the use of different weight functions for different components of the spectral density matrix. This flexibility accommodates the necessity, often observed in practice, of modeling different characteristics of the spectral density matrix components (including the cross spectrum ones) with different weight functions. Types and properties of the different weight functions are subject of most textbooks in spectral analysis and will not be discussed here. See, for instance, Priestley (1981) and references therein.

In the presence of long-range dependence, the spectral density function has a pole at the zero frequency, so that some authors restrict the summation on \eqref{smoothp} to $k\neq-j$. In practice, however, since the sample paths of $\bs X_{\!t}$ are finite with probability one, there is no problem in applying \eqref{smoothp} as defined.

Assume that the weight functions $\big(W_n^{ij}(k)\big)_{i,j=1}^q$ and the sequence $\{\ell(k)\}_{k=0}^\infty$ satisfy the following conditions:
\begin{itemize}
\item[\textbf{C}1.] $1/\ell(n)+\ell(n)/n\longrightarrow 0$, as $n$ tends to infinity;
\item[\textbf{C}2.]$W_n^{ij}(k)=W_n^{ij}(-k)$ and $W_n^{ij}(k)\geq 0$, for all $k$;
\item[\textbf{C}3.] $\sum_{|k|\leq \ell(n)}W_n^{ij}(k)=1$;
\item[\textbf{C}4.] $\sum_{|k|\leq \ell(n)}W_n^{ij}(k)^2\longrightarrow 0$, as $n$ tends to infinity.
\end{itemize}
It can be shown that under assumptions \textbf C1 - \textbf C4 and if $\bs d\in(-0.5,0)^q$, then the smoothed periodogram is a $n^{1/2}$-consistent estimator of the spectral density. Theorem \ref{consistency} thus applies and we conclude that the estimator \eqref{estimator} based on the smoothed periodogram is consistent for all $\do\in\Omega_{1/2}=(-1/4,0)^q$. At this moment, we have not been able to show the conditions of Theorems \ref{nonbeta} and \ref{norm} hold for the smoothed periodogram, but we have empirical evidence that the estimator is indeed consistent and asymptotically normally distributed for $\bs d\in(-0.5,0.5)^q$. See Section \ref{results}.

\subsection{The Tapered Periodogram}

In case of long-range dependent components in the process $\{\bs{X}_{\!t}\}_{t=0}^\infty$, the ordinary periodogram is not only non-consistent, but it is also asymptotically biased (cf. Hurvich and Beltr\~ao, 1993). A simple way to reduce this asymptotic bias is by tapering the data prior calculating the periodogram of the series. Let $\{\bs{X}_{\!t}\}_{t=0}^\infty$ be a $q$-dimensional process specified by \eqref{specif}. Let $\{h_i\}_{i=1}^q$ be a collection of real functions defined on $[0,1]$. Consider the function $L_n:\R\rightarrow\R^q$ given by $L_n(\l)=\big(L_n^1(\l), \cdots,L_n^q(\l)\big)$ where $L_n^{i}(\l):=h_{i}\big(\l/n\big)$ and let
\small\[\ S_n(\l):=\Bigg(\frac{L_n^{i}(\l)}{\sqrt{\sum_{t=1}^nL_n^{i}(t)^2}}\Bigg)_{i=1}^q.\]\normalsize
The tapered periodogram $I_T(\l;n)$ of the time series $\{\bs X_{\!t}\}_{t=1}^n$ is defined by setting
\small\begin{equation}\label{taper}
I_T(\l;n):=w_T(\l;n)\comp{w_T(\l;n)},\quad\mbox{ where }\quad w_T(\l;n):=\frac1{\sqrt{2\pi}}\sum_{t=1}^nS_n(t)\odot \bs X_{\!t}\mathrm{e}^{-it\l} .
\end{equation}\normalsize
We shall assume the following:
\begin{itemize}
\item Assumption \textbf{D}. The tapering functions $h_{i}$ are of bounded variation and $H_{i}:=\int_0^1h_{i}^2(x)dx>0$, for all $i\in\{1,\cdots,q\}.$
\end{itemize}
The tapered periodogram is not a consistent estimator of the spectral density function, since the reduction on the bias induces, in this case, an augmentation of the variance. Just as the ordinary periodogram, the increase in the variance can be dealt by smoothing the tapered periodogram in order to obtain a consistent estimator of the spectral density function in the case $\bs d\in(-0.5,0)$ (see, for instance, the recent work of Fryzlewicz et al., 2008). More details can be found in Priestley (1981), Dahlhaus (1983), Hurvich and Beltr\~ao (1993), Fryzlewicz et al. (2008) and references therein.

Under Assumption \textbf D, $\sum_{t=1}^nL_n^{i}(t)^2\sim nH_{i}$ (cf. Fryzlewicz et al., 2008) so that $I_T(\l;n)=O\big(I_n(\l)\big)$. This allows to show that the estimator  \eqref{estimator} based on the tapered periodogram is also consistent and asymptotically normally distributed. These are the contents of the next Corollaries.

        \begin{cor}\label{tapered.cons}
        Let $\{\bs{X}_{\!t}\}_{t=0}^\infty$ be a weakly stationary $q$-dimensional process specified by \eqref{specif} and with spectral density function $f$ satisfying Assumptions \textbf{A1}-\textbf{A4}. Let $f_n$ be the tapered periodogram defined in \eqref{taper} satisfying Assumption \textbf{D}. For $\do\in \Theta$, consider the estimator $\d$ based on $f_n$, as given in \eqref{estimator}. Then, $\widehat{\bs d}\overset{\P}{-\!\!\!\longrightarrow}\do$, as $n$ tends to infinity.
        \end{cor}

        \begin{cor}\label{tapered.an}
        Let $\{\bs{X}_{\!t}\}_{t=0}^\infty$ be a weakly stationary $q$-dimensional process specified by \eqref{specif} and with spectral density function $f$ satisfying Assumptions \textbf{B1}-\textbf{B5}, with \textbf{B4} holding for $\alpha=1$. Let $f_n$ be the tapered periodogram given in \eqref{taper} satisfying Assumption \textbf{D}. For $\do\in \Theta$, consider the estimator $\d$ based on $f_n$, as given in \eqref{estimator}. Then, for $\do\in\Theta$,
        \small\[m^{1/2}(\d-\do)\overset{d}{-\!\!\!\longrightarrow} N(\bs 0, \Omega),\]\normalsize
        as $n$ tends to infinity, with $\Omega$ as given in Theorem \ref{norm}.
        \end{cor}

\subsection{Simulation Results}\label{results}

In this section we present a Monte Carlo simulation study to assess the finite sample performance of the estimator \eqref{estimator}.  Recall that a $q$-dimensional stationary process $\{\bs X_t\}_{t \in \Z}$ is called a VARFIMA$(p,\bs d,q)$ if it is a stationary solution of the difference equations
\small\[\bs{\Phi}({\cal B})\,\mathrm{diag}\big\{(1-{\cal B})^{\bs d}\big\}\big(\bs X_{t}-\E(\bs X_{t})\big)=\bs{\Theta}({\cal B})\bs{\varepsilon}_{t},\]\normalsize
where $\cal B $ is the backward shift operator, $\{\bs\eps_{t}\}_{t\in\Z}$ is a $q$-dimensional stationary process (the innovation process), $\bs{\Phi}({\cal B})$ and $\bs{\Theta}({\cal B})$ are $q\times q$ matrices in ${\cal B}$, given by the equations
\small\begin{eqnarray}
\bs{\Phi}({\cal B})=\sum_{\ell=0}^{p}\bs{\phi}_{\ell}{\cal B}^\ell\,\,\,\,
\mbox{and}\,\,\,\,
\bs{\Theta}({\cal B})=\sum_{r=0}^{q}\bs{\theta}_{r}{\cal B}^r,\nonumber
\end{eqnarray}\normalsize
assumed to have no common roots, where $\bs{\phi}_{1},\cdots,\bs{\phi}_{p}, \bs{\theta}_{1},\cdots, \bs{\theta}_{q}$ are real $q\times q$ matrices and $\bs{\phi}_0 = \bs{\theta}_0 = \mathrm{I}_{q}$.

All Monte Carlo simulations are based on time series of fixed sample size $n=1,000$ obtained from bidimensional Gaussian VARFIMA$(0,\bs d,0)$ processes for several different parameters $\bs d$ and correlation $\rho\in\{0,0.3,0.6,0.8\}$. We perform 1,000 replications of each experiment. To generate the time series, we apply the traditional method of truncating  the multidimensional infinite moving average representation of the process. The truncation point is fixed in 50,000 for all cases. For comparison purposes, we calculate the estimator \eqref{shimest} (denoted by Sh) and the estimator \eqref{estimator} with the smoothed periodogram with and without the restriction $k\neq-j$ (denoted by SSh and SS$\mathrm{h}^\ast$, respectively) and with the tapered periodogram (denoted by TSh).

For the smoothed periodogram, we apply the same weights for all spectral density components, given by the so-called Bartlett's window, that is,
\small\[W_n^{ij}(k):=\frac{\sin^2(\ell(n)k/2)}{n\ell(n)\sin^2(k/2)},\quad \mbox{ for all } i,j=1,2.\]\normalsize
For the tapered periodogram, we apply the cosine-bell tapering function, namely,
\small\[h_{i}(u)=\left\{\begin{array}{cc}
                \frac12\big[1-\cos(2\pi u)\big], & \mbox{ if }\ 0\leq u\leq1/2, \vspace{.2cm}\\
                h_{i}(1-u), & \mbox{ if }\ 1/2<u\leq 1,
              \end{array}\right. \quad \mbox{ for all } i=1,2. \]\normalsize
The cosine-bell taper is widely used in applications as, for instance, in Hurvich and Ray (1995), Velasco (1999) and Olbermann et al. (2006).

The truncation point of the smoothed periodogram function is of the form $\ell(n,\beta):=\lfloor n^{\beta}\rfloor$, for $\beta\in\{0.7,0.9\}$, while the truncation point of the estimator \eqref{estimator} is of the form $m:=m(n,\alpha)=\lfloor n^{\alpha}\rfloor$, for $\alpha\in\{0.65, 0.85\}$ for all estimators. The routines are implemented in FORTRAN 95 language optimized by using OpenMP directives for parallel computing. All simulations were performed by using the computational resources from the (Brazilian) National Center of Super Computing (CESUP-UFRGS).

Tables \ref{tab+1} and \ref{tab+2} report the simulation results.  Presented are the estimated values (mean), their standard deviations (st.d.) and the mean square error of the estimates (mse). Overall all estimators present a good performance with small mse and standard deviation. The bias is generally small, except when the correlation in the noise is very high ($\rho=0.8$) and the respective component of $\bs d$ is small (specially 0.1), in which case the bias is high.  The SSh and SS$\mathrm{h}^\ast$ estimators generally perform better than Sh and TSh in terms of both, mse and bias. The same can be said about the standard deviations of the estimators. As the correlation in the innovation process increases, the estimated values degrade in some degree according to the magnitude of the respective parameter component.

The best performance in terms of bias is obtained for $\alpha=0.85$ for most cases (83 out of 128 cases). The value $\alpha=0.85$ also gives uniformly smaller mean square errors for all estimators. For the SSh estimator, there is an overall equilibrium over the values of $\beta$ presenting the smallest bias and usually the combination $\alpha=0.85$ and $\beta=0.9$ gives the best results in terms of mean square error. For the SS$\mathrm{h}^{\ast}$, $\alpha=0.85$ gives the best results in most cases (21 out of 32) while for $\beta$ there is an equilibrium between the values. The Sh and the TSh estimators present similar behavior and in most cases they agree in the value of $\alpha$ which yields smallest bias. There is a small advantage for  $\alpha=0.85$ in terms of bias in both cases (19 and 18 out of 32 for the Sh and TSh, respectively). In terms of mean square error and standard deviation, $\alpha=0.85$ generally presents best results and overall the Sh has advantage over the TSh estimator. The variance of the estimators generally responds strongly to changes in $\alpha$ than in $\beta$ in the opposite direction, that is, the higher the $\alpha$, the smaller the variance.

\begin{table}[!h]
\renewcommand{\arraystretch}{1.3}
\setlength{\tabcolsep}{3pt}
\caption{Simulation results of the estimator \eqref{estimator} based on the smoothed periodogram (SSh and SS$\mathrm h^\ast$), the ordinary periodogram (Sh) and the tapered periodogram (TSh) in VARFIMA$(0,\bs d,0)$ processes. Presented are the estimated values (mean), its standard deviation (st.d) and the mean square error of the estimates (mse).}\label{tab+1}
\vskip.3cm
\centering
{\scriptsize
\begin{tabular}{|c|c|c|c|ccc|ccc|ccc|ccc|}
\hline \hline
      \multirow{3}{*}{$\rho$}&\multirow{3}{*}{Method}&\multirow{3}{*}{$\beta$}&   \multirow{3}{*}{$\hat d_i$}    &\multicolumn{6}{|c|}{$\bs d=(0.1,0.4)$}
      &\multicolumn{6}{c|}{$\bs d=(0.2,0.3)$}   \\
      \cline{5-16}
      &       &       &      &\multicolumn{3}{|c|}{$\alpha=0.65$}       &\multicolumn{3}{|c|}{$\alpha=0.85$}            &\multicolumn{3}{|c|}{$\alpha=0.65$}
      &\multicolumn{3}{|c|}{$\alpha=0.85$}    \\
      \cline{5-16}
      &       &       &      & mean & st.d. & mse& mean & st.d. & mse& mean & st.d. & mse& mean & st.d. & mse\\
      \hline\hline
      \multirow{12}{*}{$0$}&\multirow{4}{*}{  SSh   }&\multirow{2}{*}{ 0.7 }  &$\hat d_1$&0.1043 & 0.0542 & 0.0029 &  0.0957 & 0.0268 & 0.0007 &  0.2089 & 0.0566 & 0.0033 &  0.1945 & 0.0277 & 0.0008\\
      && & $\hat d_2$ &  0.4453 & 0.0787 & 0.0082 &  0.4136 & 0.0401 & 0.0018 &  0.3160 & 0.0640 & 0.0043 &  0.2998 & 0.0320 & 0.0010\\
      \cline{3-16}
      &&\multirow{2}{*}{ 0.9   }  &$\hat d_1$&0.1062 & 0.0558 & 0.0032 &  0.0956 & 0.0267 & 0.0007 &  0.2091 & 0.0561 & 0.0032 &  0.1924 & 0.0268 & 0.0008\\
      && & $\hat d_2$ &  0.4216 & 0.0622 & 0.0043 &  0.3948 & 0.0306 & 0.0010 &  0.3074 & 0.0577 & 0.0034 &  0.2918 & 0.0288 & 0.0009\\
      \cline{2-16}
      &\multirow{4}{*}{ SS$\mathrm h^\ast$   }&\multirow{2}{*}{ 0.7 }  &$\hat d_1$&0.0940 & 0.0547 & 0.0030 &  0.0920 & 0.0270 & 0.0008 &  0.1929 & 0.0564 & 0.0032 &  0.1885 & 0.0275 & 0.0009\\
      && & $\hat d_2$ &  0.3823 & 0.0616 & 0.0041 &  0.3853 & 0.0308 & 0.0012 &  0.2783 & 0.0590 & 0.0040 &  0.2849 & 0.0298 & 0.0011\\
      \cline{3-16}
      &&\multirow{2}{*}{ 0.9   }  &$\hat d_1$&0.0995 & 0.0571 & 0.0033 &  0.0934 & 0.0270 & 0.0008 &  0.1998 & 0.0574 & 0.0033 &  0.1892 & 0.0271 & 0.0009\\
      && & $\hat d_2$ &  0.3882 & 0.0605 & 0.0038 &  0.3813 & 0.0295 & 0.0012 &  0.2866 & 0.0595 & 0.0037 &  0.2843 & 0.0292 & 0.0011\\
      \cline{2-16}
      &\multirow{2}{*}{  Sh   }& \multirow{2}{*}{  -   } &$\hat d_1$&0.1069 & 0.0577 & 0.0034 &  0.0958 & 0.0270 & 0.0007 &  0.2076 & 0.0576 & 0.0034 &  0.1916 & 0.0271 & 0.0008\\
      && & $\hat d_2$ &  0.3940 & 0.0603 & 0.0037 &  0.3820 & 0.0292 & 0.0012 &  0.2933 & 0.0598 & 0.0036 &  0.2859 & 0.0291 & 0.0010\\
      \cline{2-16}
      &\multirow{2}{*}{  TSh   }& \multirow{2}{*}{  -   } &$\hat d_1$&0.1090 & 0.0777 & 0.0061 &  0.0962 & 0.0376 & 0.0014 &  0.2087 & 0.0773 & 0.0060 &  0.1917 & 0.0376 & 0.0015\\
      && & $\hat d_2$ &  0.4051 & 0.0767 & 0.0059 &  0.3863 & 0.0384 & 0.0017 &  0.2987 & 0.0764 & 0.0058 &  0.2880 & 0.0385 & 0.0016\\
\hline \hline
      \multirow{12}{*}{$0.3$}&\multirow{4}{*}{  SSh   }&\multirow{2}{*}{ 0.7 }  &$\hat d_1$&0.1272 & 0.0526 & 0.0035 &  0.1114 & 0.0255 & 0.0008 &  0.2175 & 0.0517 & 0.0030 &  0.1998 & 0.0257 & 0.0007\\
      && & $\hat d_2$ &  0.4288 & 0.0764 & 0.0067 &  0.4054 & 0.0394 & 0.0016 &  0.3088 & 0.0587 & 0.0035 &  0.2957 & 0.0302 & 0.0009\\
      \cline{3-16}
      &&\multirow{2}{*}{ 0.9   }  &$\hat d_1$&0.1281 & 0.0537 & 0.0037 &  0.1102 & 0.0253 & 0.0007 &  0.2167 & 0.0513 & 0.0029 &  0.1972 & 0.0250 & 0.0006\\
      && & $\hat d_2$ &  0.4057 & 0.0600 & 0.0036 &  0.3864 & 0.0296 & 0.0011 &  0.3009 & 0.0530 & 0.0028 &  0.2879 & 0.0270 & 0.0009\\
      \cline{2-16}
      &\multirow{4}{*}{ SS$\mathrm h^\ast$   }&\multirow{2}{*}{ 0.7 }  &$\hat d_1$&0.1141 & 0.0523 & 0.0029 &  0.1064 & 0.0255 & 0.0007 &  0.1991 & 0.0513 & 0.0026 &  0.1931 & 0.0255 & 0.0007\\
      && & $\hat d_2$ &  0.3660 & 0.0598 & 0.0047 &  0.3769 & 0.0298 & 0.0014 &  0.2727 & 0.0544 & 0.0037 &  0.2811 & 0.0279 & 0.0011\\
      \cline{3-16}
      &&\multirow{2}{*}{ 0.9   }  &$\hat d_1$&0.1196 & 0.0544 & 0.0033 &  0.1073 & 0.0255 & 0.0007 &  0.2061 & 0.0524 & 0.0028 &  0.1937 & 0.0252 & 0.0007\\
      && & $\hat d_2$ &  0.3719 & 0.0585 & 0.0042 &  0.3727 & 0.0282 & 0.0015 &  0.2808 & 0.0550 & 0.0034 &  0.2805 & 0.0273 & 0.0011\\
      \cline{2-16}
      &\multirow{2}{*}{  Sh   }& \multirow{2}{*}{  -   } &$\hat d_1$&0.1263 & 0.0548 & 0.0037 &  0.1094 & 0.0255 & 0.0007 &  0.2135 & 0.0526 & 0.0029 &  0.1959 & 0.0252 & 0.0007\\
      && & $\hat d_2$ &  0.3781 & 0.0583 & 0.0039 &  0.3735 & 0.0279 & 0.0015 &  0.2878 & 0.0552 & 0.0032 &  0.2822 & 0.0272 & 0.0011\\
      \cline{2-16}
      &\multirow{2}{*}{  TSh   }& \multirow{2}{*}{  -   } &$\hat d_1$&0.1301 & 0.0747 & 0.0065 &  0.1110 & 0.0361 & 0.0014 &  0.2159 & 0.0712 & 0.0053 &  0.1967 & 0.0353 & 0.0013\\
      && & $\hat d_2$ &  0.3897 & 0.0735 & 0.0055 &  0.3773 & 0.0366 & 0.0019 &  0.2926 & 0.0704 & 0.0050 &  0.2838 & 0.0359 & 0.0016\\
\hline \hline
      \multirow{12}{*}{$0.6$}&\multirow{4}{*}{  SSh   }&\multirow{2}{*}{ 0.7 }  &$\hat d_1$&0.1850 & 0.0488 & 0.0096 &  0.1562 & 0.0241 & 0.0037 &  0.2353 & 0.0450 & 0.0033 &  0.2126 & 0.0226 & 0.0007\\
      && & $\hat d_2$ &  0.3987 & 0.0709 & 0.0050 &  0.3915 & 0.0383 & 0.0015 &  0.2968 & 0.0513 & 0.0026 &  0.2876 & 0.0269 & 0.0009\\
      \cline{3-16}
      &&\multirow{2}{*}{ 0.9   }  &$\hat d_1$&0.1817 & 0.0482 & 0.0090 &  0.1512 & 0.0231 & 0.0032 &  0.2323 & 0.0440 & 0.0030 &  0.2089 & 0.0219 & 0.0006\\
      && & $\hat d_2$ &  0.3780 & 0.0557 & 0.0036 &  0.3724 & 0.0279 & 0.0015 &  0.2900 & 0.0464 & 0.0023 &  0.2801 & 0.0238 & 0.0010\\
      \cline{2-16}
      &\multirow{4}{*}{ SS$\mathrm h^\ast$   }&\multirow{2}{*}{ 0.7 }  &$\hat d_1$&0.1609 & 0.0474 & 0.0060 &  0.1463 & 0.0232 & 0.0027 &  0.2118 & 0.0444 & 0.0021 &  0.2042 & 0.0223 & 0.0005\\
      && & $\hat d_2$ &  0.3391 & 0.0556 & 0.0068 &  0.3629 & 0.0282 & 0.0022 &  0.2633 & 0.0478 & 0.0036 &  0.2736 & 0.0247 & 0.0013\\
      \cline{3-16}
      &&\multirow{2}{*}{ 0.9   }  &$\hat d_1$&0.1666 & 0.0487 & 0.0068 &  0.1460 & 0.0230 & 0.0026 &  0.2190 & 0.0453 & 0.0024 &  0.2045 & 0.0221 & 0.0005\\
      && & $\hat d_2$ &  0.3448 & 0.0539 & 0.0059 &  0.3584 & 0.0260 & 0.0024 &  0.2712 & 0.0481 & 0.0031 &  0.2730 & 0.0240 & 0.0013\\
      \cline{2-16}
      &\multirow{2}{*}{  Sh   }& \multirow{2}{*}{  -   } &$\hat d_1$&0.1727 & 0.0490 & 0.0077 &  0.1475 & 0.0229 & 0.0028 &  0.2261 & 0.0454 & 0.0027 &  0.2065 & 0.0220 & 0.0005\\
      && & $\hat d_2$ &  0.3518 & 0.0537 & 0.0052 &  0.3593 & 0.0256 & 0.0023 &  0.2786 & 0.0483 & 0.0028 &  0.2748 & 0.0238 & 0.0012\\
      \cline{2-16}
      &\multirow{2}{*}{  TSh  }& \multirow{2}{*}{  -   } &$\hat d_1$&0.1802 & 0.0662 & 0.0108 &  0.1510 & 0.0327 & 0.0037 &  0.2300 & 0.0604 & 0.0045 &  0.2081 & 0.0309 & 0.0010\\
      && & $\hat d_2$ &  0.3638 & 0.0677 & 0.0059 &  0.3628 & 0.0339 & 0.0025 &  0.2831 & 0.0611 & 0.0040 &  0.2760 & 0.0316 & 0.0016\\
\hline \hline
      \multirow{12}{*}{$0.8$}&\multirow{4}{*}{  SSh   }&\multirow{2}{*}{ 0.7 }  &$\hat d_1$&0.2373 & 0.0522 & 0.0216 &  0.2099 & 0.0274 & 0.0128 &  0.2491 & 0.0445 & 0.0044 &  0.2252 & 0.0220 & 0.0011\\
      && & $\hat d_2$ &  0.3944 & 0.0697 & 0.0049 &  0.3982 & 0.0388 & 0.0015 &  0.2945 & 0.0498 & 0.0025 &  0.2858 & 0.0258 & 0.0009\\
      \cline{3-16}
      &&\multirow{2}{*}{ 0.9   }  &$\hat d_1$&0.2295 & 0.0494 & 0.0192 &  0.1998 & 0.0238 & 0.0105 &  0.2443 & 0.0424 & 0.0038 &  0.2202 & 0.0207 & 0.0008\\
      && & $\hat d_2$ &  0.3746 & 0.0561 & 0.0038 &  0.3776 & 0.0280 & 0.0013 &  0.2875 & 0.0449 & 0.0022 &  0.2780 & 0.0225 & 0.0010\\
      \cline{2-16}
      &\multirow{4}{*}{ SS$\mathrm h^\ast$   }&\multirow{2}{*}{ 0.7 }  &$\hat d_1$&0.2016 & 0.0489 & 0.0127 &  0.1937 & 0.0237 & 0.0093 &  0.2217 & 0.0433 & 0.0023 &  0.2150 & 0.0212 & 0.0007\\
      && & $\hat d_2$ &  0.3364 & 0.0562 & 0.0072 &  0.3693 & 0.0283 & 0.0017 &  0.2614 & 0.0464 & 0.0036 &  0.2718 & 0.0233 & 0.0013\\
      \cline{3-16}
      &&\multirow{2}{*}{ 0.9   }  &$\hat d_1$&0.2075 & 0.0495 & 0.0140 &  0.1915 & 0.0230 & 0.0089 &  0.2289 & 0.0439 & 0.0028 &  0.2149 & 0.0209 & 0.0007\\
      && & $\hat d_2$ &  0.3417 & 0.0541 & 0.0063 &  0.3634 & 0.0256 & 0.0020 &  0.2689 & 0.0464 & 0.0031 &  0.2709 & 0.0224 & 0.0013\\
      \cline{2-16}
      &\multirow{2}{*}{  Sh   }& \multirow{2}{*}{  -   } &$\hat d_1$&0.2137 & 0.0498 & 0.0154 &  0.1925 & 0.0229 & 0.0091 &  0.2358 & 0.0440 & 0.0032 &  0.2167 & 0.0208 & 0.0007\\
      && & $\hat d_2$ &  0.3489 & 0.0540 & 0.0055 &  0.3642 & 0.0252 & 0.0019 &  0.2764 & 0.0465 & 0.0027 &  0.2728 & 0.0223 & 0.0012\\
      \cline{2-16}
      &\multirow{2}{*}{  TSh   }& \multirow{2}{*}{  -   } &$\hat d_1$&0.2248 & 0.0659 & 0.0199 &  0.1974 & 0.0322 & 0.0105 &  0.2409 & 0.0572 & 0.0049 &  0.2187 & 0.0289 & 0.0012\\
      && & $\hat d_2$ &  0.3614 & 0.0681 & 0.0061 &  0.3681 & 0.0335 & 0.0021 &  0.2814 & 0.0587 & 0.0038 &  0.2741 & 0.0297 & 0.0016\\
\hline\hline
\end{tabular}}
\end{table}
\FloatBarrier

 \begin{table}[!h]
\renewcommand{\arraystretch}{1.3}
\setlength{\tabcolsep}{3pt}
\caption{Simulation results of the estimator \eqref{estimator} based on the smoothed periodogram (SSh and SS$\mathrm h^\ast$), the ordinary periodogram (Sh) and the tapered periodogram (TSh) in VARFIMA$(0,\bs d,0)$ processes. Presented are the estimated values (mean), its standard deviation (st.d) and the mean square error of the estimates (mse).}\label{tab+2}
\vskip.3cm
\centering
{\scriptsize
\begin{tabular}{|c|c|c|c|ccc|ccc|ccc|ccc|}
\hline \hline
      \multirow{3}{*}{$\rho$}&\multirow{3}{*}{Method}&\multirow{3}{*}{$\beta$}&   \multirow{3}{*}{$\hat d_i$}    &\multicolumn{6}{|c|}{$\bs d=(0.1,0.3)$}
      &\multicolumn{6}{c|}{$\bs d=(0.3,0.4)$}   \\
      \cline{5-16}
      &       &       &      &\multicolumn{3}{|c|}{$\alpha=0.65$}       &\multicolumn{3}{|c|}{$\alpha=0.85$}            &\multicolumn{3}{|c|}{$\alpha=0.65$}
      &\multicolumn{3}{|c|}{$\alpha=0.85$}    \\
      \cline{5-16}
      &       &       &      & mean & st.d. & mse& mean & st.d. & mse& mean & st.d. & mse& mean & st.d. & mse\\
      \hline\hline
      \multirow{12}{*}{$0$}&\multirow{4}{*}{  SSh   }&\multirow{2}{*}{ 0.7 }  &$\hat d_1$&0.1038 & 0.0541 & 0.0029 &  0.0956 & 0.0268 & 0.0007 &  0.3237 & 0.0610 & 0.0043 &  0.2987 & 0.0299 & 0.0009\\
      && & $\hat d_2$ &  0.3156 & 0.0641 & 0.0044 &  0.2998 & 0.0321 & 0.0010 &  0.4454 & 0.0782 & 0.0082 &  0.4136 & 0.0401 & 0.0018\\
      \cline{3-16}
      &&\multirow{2}{*}{ 0.9   }  &$\hat d_1$&0.1057 & 0.0558 & 0.0031 &  0.0955 & 0.0267 & 0.0007 &  0.3173 & 0.0570 & 0.0035 &  0.2916 & 0.0273 & 0.0008\\
      && & $\hat d_2$ &  0.3068 & 0.0577 & 0.0034 &  0.2918 & 0.0288 & 0.0009 &  0.4225 & 0.0620 & 0.0044 &  0.3948 & 0.0306 & 0.0010\\
      \cline{2-16}
      &\multirow{4}{*}{ SS$\mathrm h^\ast$   }&\multirow{2}{*}{ 0.7 }  &$\hat d_1$&0.0935 & 0.0546 & 0.0030 &  0.0919 & 0.0270 & 0.0008 &  0.2975 & 0.0586 & 0.0034 &  0.2880 & 0.0284 & 0.0010\\
      && & $\hat d_2$ &  0.2777 & 0.0591 & 0.0040 &  0.2848 & 0.0298 & 0.0011 &  0.3831 & 0.0613 & 0.0040 &  0.3854 & 0.0308 & 0.0012\\
      \cline{3-16}
      &&\multirow{2}{*}{ 0.9   }  &$\hat d_1$&0.0990 & 0.0570 & 0.0032 &  0.0933 & 0.0270 & 0.0008 &  0.3032 & 0.0581 & 0.0034 &  0.2863 & 0.0274 & 0.0009\\
      && & $\hat d_2$ &  0.2859 & 0.0596 & 0.0037 &  0.2843 & 0.0292 & 0.0011 &  0.3893 & 0.0603 & 0.0037 &  0.3814 & 0.0294 & 0.0012\\
      \cline{2-16}
      &\multirow{2}{*}{  Sh   }& \multirow{2}{*}{  -   } &$\hat d_1$&0.1063 & 0.0576 & 0.0034 &  0.0957 & 0.0270 & 0.0007 &  0.3103 & 0.0579 & 0.0035 &  0.2881 & 0.0272 & 0.0009\\
      && & $\hat d_2$ &  0.2926 & 0.0599 & 0.0036 &  0.2859 & 0.0291 & 0.0010 &  0.3951 & 0.0601 & 0.0036 &  0.3821 & 0.0292 & 0.0012\\
      \cline{2-16}
      &\multirow{2}{*}{  TSh   }& \multirow{2}{*}{  -   } &$\hat d_1$&0.1079 & 0.0775 & 0.0061 &  0.0960 & 0.0376 & 0.0014 &  0.3124 & 0.0770 & 0.0061 &  0.2884 & 0.0375 & 0.0015\\
      && & $\hat d_2$ &  0.2976 & 0.0765 & 0.0059 &  0.2879 & 0.0385 & 0.0016 &  0.4071 & 0.0765 & 0.0059 &  0.3865 & 0.0384 & 0.0017\\
\hline \hline
      \multirow{12}{*}{$0.3$}&\multirow{4}{*}{  SSh   }&\multirow{2}{*}{ 0.7 }  &$\hat d_1$&0.1204 & 0.0506 & 0.0030 &  0.1061 & 0.0251 & 0.0007 &  0.3333 & 0.0557 & 0.0042 &  0.3043 & 0.0277 & 0.0008\\
      && & $\hat d_2$ &  0.3022 & 0.0603 & 0.0036 &  0.2927 & 0.0307 & 0.0010 &  0.4380 & 0.0719 & 0.0066 &  0.4096 & 0.0381 & 0.0015\\
      \cline{3-16}
      &&\multirow{2}{*}{ 0.9   }  &$\hat d_1$&0.1213 & 0.0520 & 0.0032 &  0.1054 & 0.0250 & 0.0007 &  0.3253 & 0.0521 & 0.0033 &  0.2965 & 0.0254 & 0.0007\\
      && & $\hat d_2$ &  0.2941 & 0.0542 & 0.0030 &  0.2848 & 0.0273 & 0.0010 &  0.4161 & 0.0569 & 0.0035 &  0.3910 & 0.0288 & 0.0009\\
      \cline{2-16}
      &\multirow{4}{*}{ SS$\mathrm h^\ast$  }&\multirow{2}{*}{ 0.7 }  &$\hat d_1$&0.1075 & 0.0508 & 0.0026 &  0.1016 & 0.0252 & 0.0006 &  0.3036 & 0.0531 & 0.0028 &  0.2923 & 0.0263 & 0.0008\\
      && & $\hat d_2$ &  0.2654 & 0.0556 & 0.0043 &  0.2779 & 0.0282 & 0.0013 &  0.3778 & 0.0565 & 0.0037 &  0.3818 & 0.0290 & 0.0012\\
      \cline{3-16}
      &&\multirow{2}{*}{ 0.9   }  &$\hat d_1$&0.1132 & 0.0530 & 0.0030 &  0.1028 & 0.0252 & 0.0006 &  0.3093 & 0.0529 & 0.0029 &  0.2906 & 0.0255 & 0.0007\\
      && & $\hat d_2$ &  0.2735 & 0.0561 & 0.0038 &  0.2773 & 0.0276 & 0.0013 &  0.3838 & 0.0556 & 0.0034 &  0.3778 & 0.0276 & 0.0013\\
      \cline{2-16}
      &\multirow{2}{*}{  Sh   }& \multirow{2}{*}{  -   } &$\hat d_1$&0.1201 & 0.0534 & 0.0033 &  0.1049 & 0.0252 & 0.0007 &  0.3161 & 0.0528 & 0.0030 &  0.2922 & 0.0253 & 0.0007\\
      && & $\hat d_2$ &  0.2806 & 0.0564 & 0.0036 &  0.2790 & 0.0274 & 0.0012 &  0.3899 & 0.0555 & 0.0032 &  0.3786 & 0.0273 & 0.0012\\
      \cline{2-16}
      &\multirow{2}{*}{  TSh   }& \multirow{2}{*}{  -   } &$\hat d_1$&0.1228 & 0.0726 & 0.0058 &  0.1061 & 0.0356 & 0.0013 &  0.3200 & 0.0711 & 0.0054 &  0.2934 & 0.0352 & 0.0013\\
      && & $\hat d_2$ &  0.2854 & 0.0717 & 0.0053 &  0.2805 & 0.0362 & 0.0017 &  0.4011 & 0.0704 & 0.0050 &  0.3824 & 0.0359 & 0.0016\\
\hline \hline
      \multirow{12}{*}{$0.6$}&\multirow{4}{*}{  SSh   }&\multirow{2}{*}{ 0.7 }  &$\hat d_1$&0.1571 & 0.0451 & 0.0053 &  0.1334 & 0.0226 & 0.0016 &  0.3543 & 0.0493 & 0.0054 &  0.3186 & 0.0247 & 0.0010\\
      && & $\hat d_2$ &  0.2806 & 0.0542 & 0.0033 &  0.2803 & 0.0282 & 0.0012 &  0.4249 & 0.0627 & 0.0046 &  0.4015 & 0.0347 & 0.0012\\
      \cline{3-16}
      &&\multirow{2}{*}{ 0.9   }  &$\hat d_1$&0.1553 & 0.0452 & 0.0051 &  0.1310 & 0.0221 & 0.0014 &  0.3425 & 0.0449 & 0.0038 &  0.3086 & 0.0223 & 0.0006\\
      && & $\hat d_2$ &  0.2736 & 0.0486 & 0.0031 &  0.2725 & 0.0246 & 0.0014 &  0.4049 & 0.0500 & 0.0025 &  0.3834 & 0.0258 & 0.0009\\
      \cline{2-16}
      &\multirow{4}{*}{  SS$\mathrm h^\ast$  }&\multirow{2}{*}{ 0.7 }  &$\hat d_1$&0.1374 & 0.0449 & 0.0034 &  0.1266 & 0.0223 & 0.0012 &  0.3163 & 0.0457 & 0.0024 &  0.3034 & 0.0229 & 0.0005\\
      && & $\hat d_2$ &  0.2462 & 0.0501 & 0.0054 &  0.2657 & 0.0256 & 0.0018 &  0.3688 & 0.0498 & 0.0035 &  0.3747 & 0.0259 & 0.0013\\
      \cline{3-16}
      &&\multirow{2}{*}{ 0.9   }  &$\hat d_1$&0.1435 & 0.0464 & 0.0040 &  0.1273 & 0.0223 & 0.0012 &  0.3220 & 0.0456 & 0.0026 &  0.3013 & 0.0222 & 0.0005\\
      && & $\hat d_2$ &  0.2538 & 0.0503 & 0.0047 &  0.2650 & 0.0247 & 0.0018 &  0.3745 & 0.0487 & 0.0030 &  0.3705 & 0.0243 & 0.0015\\
      \cline{2-16}
      &\multirow{2}{*}{  Sh   }& \multirow{2}{*}{  -   } &$\hat d_1$&0.1500 & 0.0468 & 0.0047 &  0.1290 & 0.0223 & 0.0013 &  0.3284 & 0.0455 & 0.0029 &  0.3026 & 0.0221 & 0.0005\\
      && & $\hat d_2$ &  0.2614 & 0.0506 & 0.0040 &  0.2669 & 0.0246 & 0.0017 &  0.3809 & 0.0486 & 0.0027 &  0.3714 & 0.0240 & 0.0014\\
      \cline{2-16}
      &\multirow{2}{*}{  TSh   }& \multirow{2}{*}{  -   } &$\hat d_1$&0.1544 & 0.0625 & 0.0069 &  0.1311 & 0.0315 & 0.0020 &  0.3353 & 0.0603 & 0.0049 &  0.3050 & 0.0309 & 0.0010\\
      && & $\hat d_2$ &  0.2663 & 0.0639 & 0.0052 &  0.2680 & 0.0325 & 0.0021 &  0.3914 & 0.0614 & 0.0038 &  0.3747 & 0.0317 & 0.0016\\
\hline \hline
      \multirow{12}{*}{$0.8$}&\multirow{4}{*}{  \emph{SSh}   }&\multirow{2}{*}{ 0.7 }  &$\hat d_1$&0.1882 & 0.0468 & 0.0100 &  0.1641 & 0.0231 & 0.0046 &  0.3710 & 0.0504 & 0.0076 &  0.3332 & 0.0253 & 0.0017\\
      && & $\hat d_2$ &  0.2787 & 0.0544 & 0.0034 &  0.2822 & 0.0280 & 0.0011 &  0.4228 & 0.0605 & 0.0042 &  0.4003 & 0.0338 & 0.0011\\
      \cline{3-16}
      &&\multirow{2}{*}{ 0.9   }  &$\hat d_1$&0.1840 & 0.0452 & 0.0091 &  0.1593 & 0.0217 & 0.0040 &  0.3562 & 0.0440 & 0.0051 &  0.3208 & 0.0215 & 0.0009\\
      && & $\hat d_2$ &  0.2713 & 0.0486 & 0.0032 &  0.2736 & 0.0239 & 0.0013 &  0.4028 & 0.0485 & 0.0024 &  0.3817 & 0.0247 & 0.0009\\
      \cline{2-16}
      &\multirow{4}{*}{  SS$\mathrm h^\ast$   }&\multirow{2}{*}{ 0.7 }  &$\hat d_1$&0.1625 & 0.0455 & 0.0060 &  0.1545 & 0.0220 & 0.0035 &  0.3263 & 0.0447 & 0.0027 &  0.3145 & 0.0218 & 0.0007\\
      && & $\hat d_2$ &  0.2443 & 0.0501 & 0.0056 &  0.2672 & 0.0249 & 0.0017 &  0.3671 & 0.0485 & 0.0034 &  0.3733 & 0.0248 & 0.0013\\
      \cline{3-16}
      &&\multirow{2}{*}{ 0.9   }  &$\hat d_1$&0.1687 & 0.0465 & 0.0069 &  0.1543 & 0.0217 & 0.0034 &  0.3318 & 0.0442 & 0.0030 &  0.3117 & 0.0210 & 0.0006\\
      && & $\hat d_2$ &  0.2514 & 0.0499 & 0.0048 &  0.2659 & 0.0237 & 0.0017 &  0.3724 & 0.0471 & 0.0030 &  0.3687 & 0.0229 & 0.0015\\
      \cline{2-16}
      &\multirow{2}{*}{  Sh   }& \multirow{2}{*}{  -   } &$\hat d_1$&0.1753 & 0.0468 & 0.0079 &  0.1558 & 0.0216 & 0.0036 &  0.3381 & 0.0440 & 0.0034 &  0.3128 & 0.0208 & 0.0006\\
      && & $\hat d_2$ &  0.2592 & 0.0502 & 0.0042 &  0.2678 & 0.0236 & 0.0016 &  0.3790 & 0.0469 & 0.0026 &  0.3695 & 0.0226 & 0.0014\\
      \cline{2-16}
      &\multirow{2}{*}{  TSh   }& \multirow{2}{*}{  -   } &$\hat d_1$&0.1811 & 0.0611 & 0.0103 &  0.1584 & 0.0303 & 0.0043 &  0.3477 & 0.0573 & 0.0055 &  0.3161 & 0.0289 & 0.0011\\
      && & $\hat d_2$ &  0.2647 & 0.0633 & 0.0053 &  0.2692 & 0.0315 & 0.0019 &  0.3901 & 0.0592 & 0.0036 &  0.3730 & 0.0299 & 0.0016\\
\hline\hline
\end{tabular}}
\end{table}
\FloatBarrier

%

\begin{figure}[h]
\centering
\mbox{
 \subfigure[]{\includegraphics[width=0.2\textwidth]{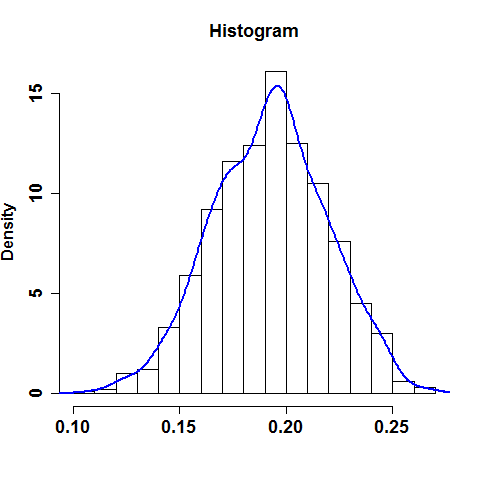}}\hskip.2cm
 \subfigure[]{\includegraphics[width=0.2\textwidth]{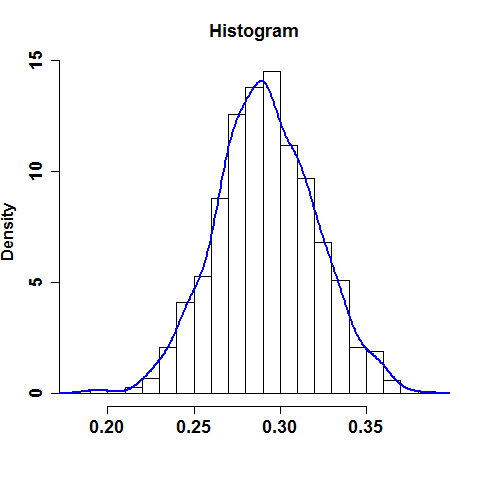}}\hskip.2cm
 \subfigure[]{\includegraphics[width=0.2\textwidth]{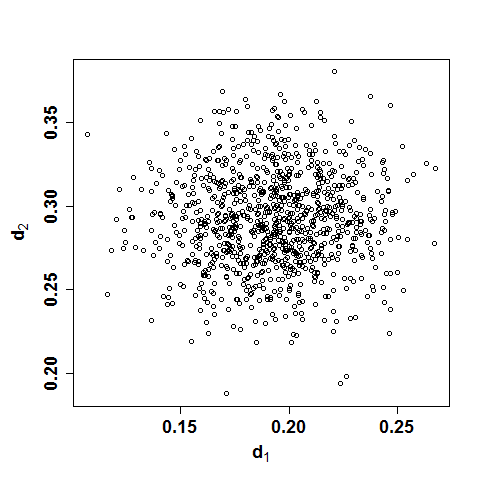}}
 }
 \mbox{
 \subfigure[]{\includegraphics[width=0.2\textwidth]{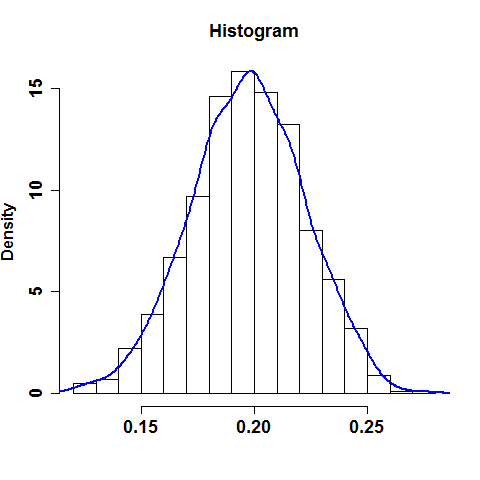}}\hskip.2cm
 \subfigure[]{\includegraphics[width=0.2\textwidth]{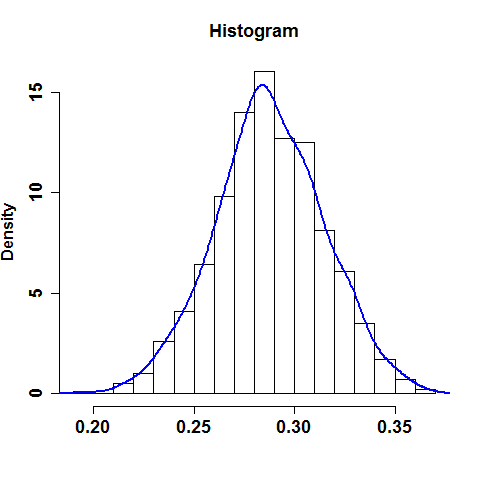}}\hskip.2cm
 \subfigure[]{\includegraphics[width=0.2\textwidth]{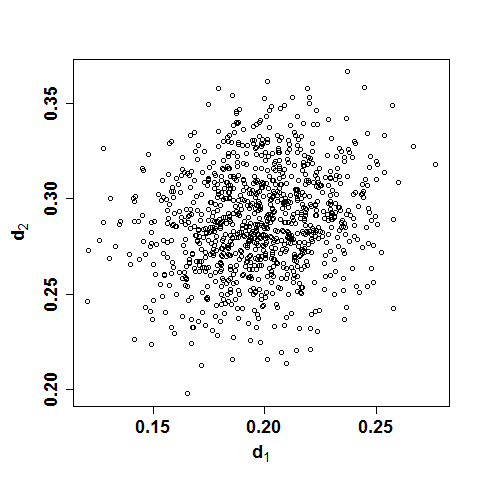}}
  }
  \mbox{
 \subfigure[]{\includegraphics[width=0.2\textwidth]{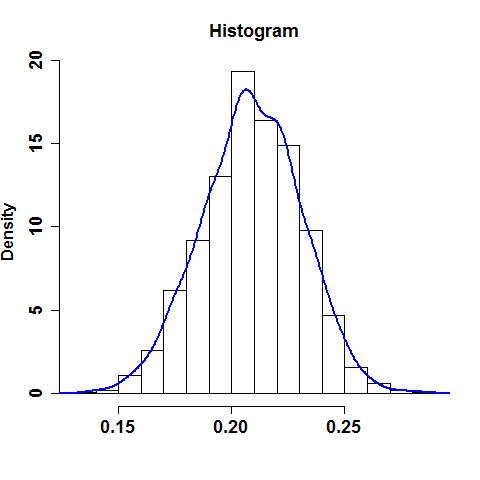}}\hskip.2cm
 \subfigure[]{\includegraphics[width=0.2\textwidth]{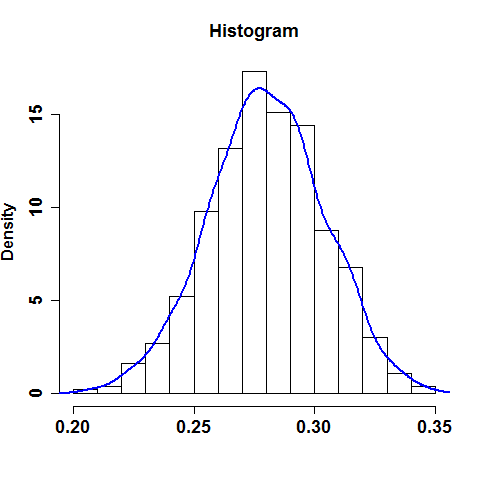}}\hskip.2cm
  \subfigure[]{\includegraphics[width=0.2\textwidth]{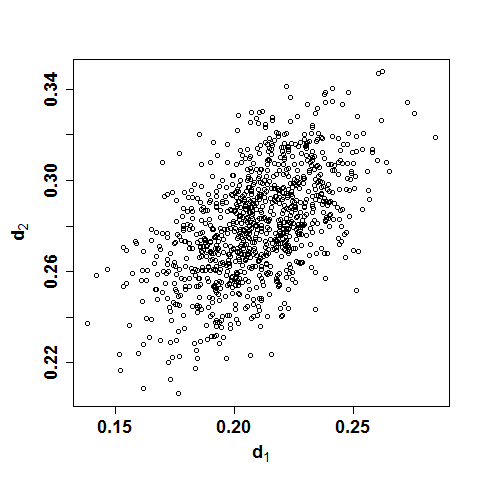}}
  }
  \mbox{
 \subfigure[]{\includegraphics[width=0.2\textwidth]{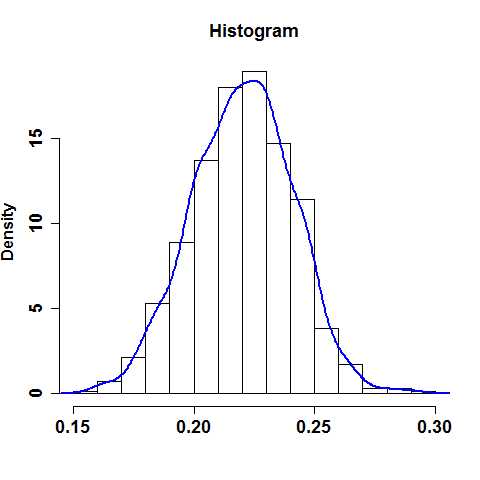}}\hskip.2cm
 \subfigure[]{\includegraphics[width=0.2\textwidth]{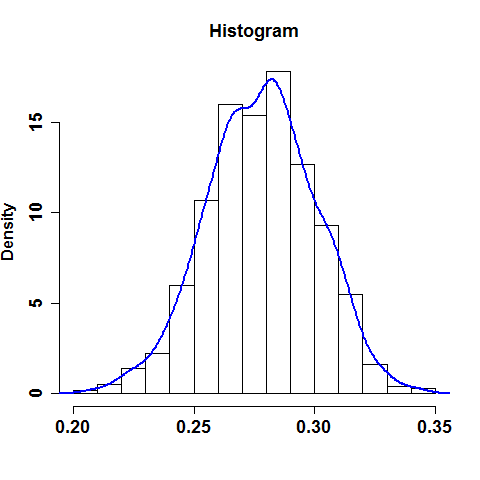}}\hskip.2cm
   \subfigure[]{\includegraphics[width=0.2\textwidth]{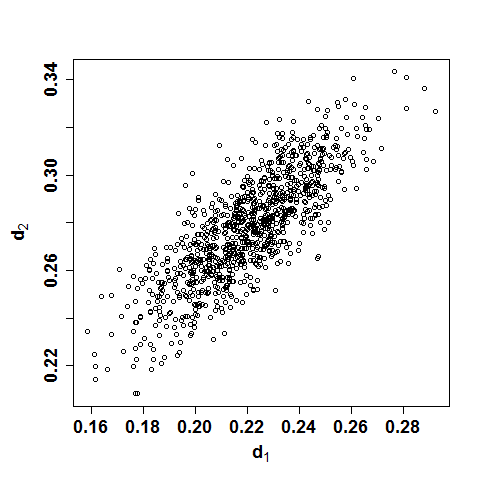}}
  }
  \caption{Histogram, kernel density and scatter plot of the SSh estimated values of $\bs d_0=(0.2,0.3)$ for (a)--(c) $\rho =0$; (d)--(f) $\rho =0.3$; (g)--(i) $\rho =0.6$ and (j)--(l) $\rho=0.8$. }\label{hist}
\end{figure}

Figure \ref{hist} presents the scatter plot, histogram and kernel density estimator of the SSh estimated values for $\bs d_0=(0.2,0.3)$ when $\alpha=0.85$ and $\beta=0.9$. Figures \ref{hist}(a)--(c) correspond to $\rho=0$, Figures \ref{hist}(d)--(f) to $\rho=0.3$, Figures \ref{hist}(g)--(i) to $\rho=0.6$ and Figures \ref{hist}(j)--(l) to $\rho=0.8$. At this moment, we were not able to prove the asymptotic normality of the SSh estimator by direct verification of \eqref{cond_an}. However, we conjecture that this is the case and Figure \ref{hist} supports this opinion.

\FloatBarrier

\section{Conclusions}

In this work we propose and analyze a class of Gaussian semiparametric estimators of multivariate long-range dependent processes. The work is motivated by the semiparametric methodology presented in Shimotsu (2007). More specifically, we propose a class of estimators based on the method studied in Shimotsu (2007) by substituting the periodogram applied there for an arbitrary spectral density estimator. We analyze two frameworks. First we assume that the spectral density estimator is consistent for the spectral density estimator and we show that the proposed semiparametric estimator is also consistent under mild conditions. Second, we relax the consistency condition and derive necessary conditions for the consistency and asymptotic normality of the proposed estimator.  We show that the variance-covariance matrix of the limiting distribution is the same as the one derived
in Shimotsu (2007), under the same conditions imposed in the process.

In order to assess the finite sample performance and illustrate the usefulness of the estimator, we perform a Monte Carlo simulation based on VARFIMA$(0,\bs d,0)$ processes. We applied the smoothed periodogram with the Bartlett's weight function and the tapered periodogram with the cosine-bell taper as the spectral density estimators. For comparison we also compute the estimator proposed in Shimotsu (2007).

The assumptions required in the asymptotic theory are mild ones and are commonly applied in the literature. The semiparametric methodology present several advantages compared to the parametric framework such as weaker distributional  assumptions, robustness with respect to misspecification of the short run dynamics of the process and efficiency. The theory includes the fractionally integrated processes as well as the class of VARFIMA processes.

\vskip.3cm

\noindent \small\textbf{Acknowledgements}\vspace{.3cm}\\
G. Pumi was partially supported by CAPES/Fulbright Grant BEX 2910/06-3 and by CNPq-Brazil. S.R.C. Lopes research was partially supported by CNPq-Brazil, by CAPES-Brazil, by Pronex {\it Probabilidade e Processos Estoc\'asticos} - E-26/170.008/2008 -APQ1 and also by INCT {\it em Matem\'atica}. The authors are also grateful to the (Brazilian) National Center of Super Computing (CESUP-UFRGS) for the computational resources. \normalsize

\section*{Appendix A: Proofs}
\renewcommand{\theequation}{A.\arabic{equation}}
\renewcommand{\thesubsection}{A}
In this section we present the proofs of the results in Sections 4 and 5. We establish lemmas and  theorems in the same sequence as they appear in the text.

\subsubsection*{Proof of Lemma \ref{lema1}:}
 By hypothesis, $f_n(\l)=f(\l)+o_\P(n^{-\beta})$ in $B$. Recalling the definition of $\L_j$ given in \eqref{fapprox}, we have
\small\begin{align}\label{lema1e1}
\widehat G(\do)&=\frac{1}{m}\sum_{j=1}^m\re\big[\L_j(\do)^{-1}f_n(\l_j)\comp{\L_j(\do)^{-1}}\big]=\frac{1}{m}\sum_{j=1}^m\re\big[\L_j(\do)^{-1}\big(f(\l_j) +o_\P(n^{-\beta})\big)\comp{\L_j(\do)^{-1}}\big]\nonumber\\
&=G_0+\frac{1}{m}\sum_{j=1}^m\re\big[\L_j(\do)^{-1}o_\P(n^{-\beta})\comp{\L_j(\do)^{-1}}\big].
\end{align}\normalsize
The $(r,s)$-th component of the second part on the RHS of \eqref{lema1e1} is given by
\small\begin{align*}
\frac1m\sum_{j=1}^m\re\big[\mathrm e^{\im(\l_j-\pi)(d_r^0-d_s^0)/2}&\l_j^{d_r^0+d_s^0}\big]o_\P(n^{-\beta})=\bigg[\frac{1}{m}\sum_{j=1}^m\l_j^{2d_r^0+d_s^0}\bigg]o_\P(n^{-\beta})\\
&=\frac{1}{{d_r^0+d_s^0}+1}\left(\frac{2\pi m}{n}\right)^{d_r^0+d_s^0}\bigg[\frac{{d_r^0+d_s^0}+1}{m}\sum_{j=1}^m\left(\frac{j}{m}\right)^{d_r^0+d_s^0}\bigg]o_\P(n^{-\beta})\\
&=\frac{1}{d_r^0+d_s^0+1}\left(\frac{2\pi m}{n}\right)^{d_r^0+d_s^0}\big[O(m^{\beta-1})+1\big]o_\P(n^{-\beta})\,\,=\,\,o_\P(1),
\end{align*}\normalsize
\pagebreak
where the penultimate equality follows from lemma 1 in Robinson (1995b), by taking $\gamma=d_r^0+d_s^0+1>1-\beta>0$, while the last one follows since $\beta\in(0,1)$  and $\do\in\Omega_\beta$. Hence, $\widehat G(\do)_{rs}=G_0^{rs}+o_\P(1)$, for all $r,s\in\{1,\cdots,q\}$ and the proof is complete.\fim

\subsubsection*{Proof of Theorem \ref{consistency}:}

Let $\bs\te=(\te_1,\cdots,\te_q)^\prime:=\bs d-\do$ and $L(\bs d):=S(\bs d)-S(\do)$.  Let $0<\delta<1/2$ be fixed and let
\small\[N_\delta:=\big\{\bs d:\|\bs d-\do\|_{\infty}>\delta\big\}.\]\normalsize
Let $0<\epsilon<1/4$ and define $\Theta_1:=\big\{\bs\te:\bs\te\in[-1/2+\epsilon,1/2]^q\big\}$ and $\Theta_2=\Omega_\beta\!\setminus\!\Theta_1$ (possibly an empty set), where $\Omega_\beta$ is given by \eqref{set}. Following Robinson (1995b) and Shimotsu (2007), we have
\small\begin{align}\label{diffp}
\P\big(\|\d-\do\|_\infty>\delta\big)&\leq\P\Big(\inf_{\overline{N_\delta}\cap\Omega_\beta}\!\!\big\{L(\bs d)\big\}\leq 0\Big)\nonumber\\
&\leq\P\Big(\inf_{\overline{N_\delta}\cap\Theta_1}\big\{L(\bs d)\big\}\leq 0\Big)+\P\Big(\inf_{\Theta_2}\big\{L(\bs d)\big\}\leq 0\Big):=P_1+P_2
\end{align}\normalsize
where, for a given set $\mathcal{O}$, $\overline{\mathcal O}$ denotes the closure of $\mathcal O$. We shall  first show that $P_1\rightarrow 0$, as $n$ tends to infinity. Rewrite $L(\bs d)$ as
\small\begin{align}\label{aux1}
L(\bs d )&=\log\big(\det \{\widehat G(\bs d)\}\big)-\log\big(\det \{\widehat G(\do)\}\big)-2\sum_{k=1}^q\te_k\frac{1}{m}\sum_{j=1}^m\log(\l_j)\nonumber\\
&=\log\big(\det \{\widehat G(\bs d)\}\big)-\log\big(\det \{\widehat G(\do)\}\big)+\log\left(\frac{2\pi m}{n}\right)^{\!\!-2\sum_{k}\te_k}-\nonumber\\
&\hspace{1cm}-2\sum_{k=1}^q\te_k\bigg(\frac{1}{m}\sum_{j=1}^m\log(j)-\log(m)\bigg)-\sum_{k=1}^q\log(2\te_k+1)\nonumber\\
&=\log\big(\mz{A}(\bs d)\big)-\log\big(\mz{B}(\bs d)\big)-\log\big(\mz{A}(\do)\big)+\log\big(\mz{B}(\do)\big)+\mz{R}(\bs d)\nonumber\\
&= Q_1(\bs d)-Q_2(\bs d)+\mz{R}(\bs d),
\end{align}\normalsize
where
\small\begin{gather*}
Q_1(\bs d):=\log\big(\mz{A}(\bs d)\big)-\log\big(\mz{B}(\bs d)\big),\qquad Q_2(\bs d):=\log\big(\mz{A}(\do)\big)+\log\big(\mz{B}(\do)\big),\\
\mz{A}(\bs d): =\left(\frac{2\pi m}{n}\right)^{\!\!-2\sum_{k}\te_k}\!\!\det \{\widehat G(\bs d)\},\qquad
\mz{B}(\bs d): =\det \{G_0\}\prod_{k=1}^q\frac{1}{2\te_k+1},\\
\mbox{ and }\quad\mz{R}(\bs d): =2\sum_{k=1}^q\te_k\bigg(\log(m)-\frac{1}{m}\sum_{j=1}^m\log(j)\bigg)-\sum_{k=1}^q\log(2\te_k+1).\phantom{ and and a}
\end{gather*}
By lemma 2 in Robinson (1995b),  $\log(m)-m^{-1}\sum_{j=1}^m\log(j)=1+O(m^{-1}\log(m))$, so that
\small\[\mz{R}(\bs d)=\sum_{k=1}^q2\te_k-\log(2\te_k+1)+O\left(\frac{\log(m)}{m}\right).\]\normalsize
Since $x-\log(x+1)$ has a unique global minimum in $(-1,\infty)$ at $x=0$ and $x-\log(x+1)\geq x^2/4$, for $|x|\leq1$, it follows that
\small\[\inf_{\overline{N_\delta}\cap\Theta_1}\!\!\big\{\mz{R}(\bs d)\big\}\geq \frac{1}{4}\Big(2\max_{k}\{\te_k\}\Big)^2\geq\delta^2>0\,.\]\normalsize
As for $Q_1(\bs d)$ and $Q_2(\bs d)$ in \eqref{aux1}, it suffices to show the existence of a function $h(\bs d)>0$ satisfying
\small\[(\mathrm{i})\,\,\sup_{\Theta_1}\big\{\big|\mz{A}(\bs d)-h(\bs d)\big|\big\}=o_{\P}(1);\quad(\mathrm{ii})\,\,h(\bs d)\geq \mz{B}(\bs d);\quad (\mathrm{iii})\,\,h(\do)=\mz{B}(\do),\]\normalsize
as $n$ goes to infinity, because (ii) implies $\displaystyle{\inf_{\Theta_1}}\big\{h(\bs d)\big\}\geq \displaystyle{\inf_{\Theta_1}}\big\{\mz{B}(\bs d)\big\}>0$, so that, uniformly in $\Theta_1$,
\small\begin{equation}\label{lcond}
Q_1(\bs d)\geq \log\big(\mz{A}(\bs d)\big)-\log\big(h(\bs d)\big)=\log\big(h(\bs d)+o_\P(1)\big)-\log\big(h(\bs d)\big)=o_\P(1),
\end{equation}\normalsize
and (iii) implies $Q_2(\bs d)=\log\big(h(\do)+o_\P(1)\big)-\log\big(h(\do)\big)=o_\P(1)$, from which $P_1\rightarrow 0$ follows. To show (i), recall that
\small\begin{align}\label{4.15A}
\L_j(\bs d)^{-1}&=\!\!\diag_{k\in\{1,\cdots,q\}}\!\!\{\l_j^{d_k}\mathrm{e}^{\im(\l_j-\pi)d_k/2}\}=\!\!
\diag_{k\in\{1,\cdots,q\}}\!\!\{\l_j^{(d_k-d_k^0)}\mathrm{e}^{\im(\l_j-\pi)(d_k-d_k^0)/2}\times\l_j^{d_k^0}\mathrm{e}^{\im(\l_j-\pi)d_k^0/2}\}\nonumber\\
&=\L_j(\bs d-\do)^{-1}\L_j(\do)^{-1}=\L_j(\bs \te)^{-1}\L_j(\do)^{-1}\!,
\end{align}\normalsize
so that we can write
\small\begin{align}\label{Mjs}
\mz A&(\bs d)=\left(\frac{2\pi m}{n}\right)^{\!\!-2\sum_{k}\te_k}\!\!\times\det \bigg\{\frac{1}{m}\sum_{j=1}^m\re\big[\L_j(\bs\te)^{-1}\L_j(\do)^{-1}f_n(\l_j)\comp{\L_j(\do)^{-1}} \comp{\L_j(\bs\te)^{-1}}\big]\bigg\}\nonumber\\
&=\det \bigg\{\frac{1}{m}\sum_{j=1}^m\re\big[M_j(\bs\te)\L_j(\do)^{-1}\big(f(\l_j)+o_\P(n^{-\beta})\big)\comp{\L_j(\do)^{-1}} \comp{M_j(\bs\te)}\big]\bigg\}\nonumber\\
&=\det \bigg\{\frac{1}{m}\sum_{j=1}^m\re\big[M_j(\bs\te)G_0\comp{M_j(\bs\te)}\big] +\frac{1}{m}\sum_{j=1}^m\re\big[M_j(\bs\te)\L_j(\do)^{-1}o_\P(n^{-\beta})\comp{\L_j(\do)^{-1}}\comp{M_j(\bs\te)}\big]\bigg\},
\end{align}\normalsize
where $M_j(\bs\te):=\diag_{k\in\{1,\cdots,q\}}\Big\{\mathrm{e}^{\im(\l_j-\pi)\te_k/2}(j/m)^{\te_k}\Big\}$. To determine the function $h$ in (i), we first show that the second term on the RHS of \eqref{Mjs} is $o_\P(1)$. This follows by noticing that its $(r,s)$-th element is given by
\small
\begin{align*}
\frac1m\sum_{j=1}^m\re\bigg[\mathrm{e}^{\im(\l_j-\pi)(\hat d_r-\hat d_s)/2}&\left(\frac jm\right)^{\te_r+\te_s}\l_j^{d_r^0+d_s^0}\bigg]\leq \frac{\mz C}m\sum_{j=1}^m\left(\frac jm\right)^{\hat d_r+\hat d_s}\bigg(\frac mn\bigg)^{d_r^0+d_s^0}o_\P(n^{-\beta})\nonumber\\
&=\frac1{\hat d_r+\hat d_s+1}\bigg(\frac{2\pi m}n\bigg)^{\hat d_r+\hat d_s}\Big[O(m^{\beta-1})+1\Big]\bigg(\frac mn\bigg)^{d_r^0+d_s^0}o_\P(n^{-\beta})=o_\P(1)
\end{align*}\normalsize
for $\mz C$ a constant, where the penultimate equality follows from lemma 1 in Robinson (1995b) and the last one follows since $\d\in\Omega_\beta$. Hence
\small\begin{align}\label{auxl4}
\mz A&(\bs d)=\det \bigg\{\frac{1}{m}\sum_{j=1}^m\re\big[M_j(\bs\te)G_0\comp{M_j(\bs\te)}\big] +o_\P(1)\bigg\}.
\end{align}\normalsize
Upon defining the matrices
\small\[\mathcal{E}(\bs \te):=\big(\mathrm{e}^{-\im\pi(\te_r-\te_s)/2}\big)_{r,s=1}^q \quad \mbox{ and }\quad \mathcal{M}(\bs\te):=\left(\frac{1}{1+\te_r+\te_s}\right)_{r,s=1}^q\]\normalsize
from the proof of theorem 1 in Shimotsu (2007), it follows that the function
\small\[h(\bs d):=\det \Big\{\re\big[\mathcal E(\bs\te)\big]\odot\mathcal M(\bs\te)\odot G_0\Big\},\]\normalsize
where $\odot$ denotes the Hadamard product, satisfies the conditions (i), (ii) and (iii) in \eqref{lcond} (see the argument following (11) in Shimotsu, 2007, p.292).

Now we move to bound $P_2$ in \eqref{diffp}. Expression \eqref{4.15A} can be used to rewrite $L(\bs d)$ as
\small\begin{eqnarray}\label{aux6}
L(\bs d )&=&\log\big(\det \{\widehat G(\bs d)\}\big)-\log\big(\det \{\widehat G(\do)\}\big)-2\sum_{k=1}^q\te_k\frac{1}{m}\sum_{j=1}^m\log(\l_j)\nonumber\\
&=&\log\Big(\det \big\{\mcd(\bs d)\big\}\Big)-\log\Big(\det \big\{\mcd(\do)\big\}\Big),
\end{eqnarray}\normalsize
where
\small\[\mcd(\bs d):=\frac{1}{m}\sum_{j=1}^m\re\big[\mcp_j(\bs \te)\L_j(\do)^{-1}f_n(\l_j)\comp{\L_j(\do)^{-1}}\comp{\mcp_j(\bs \te)}\big],\]\normalsize
with
\small\[  \mcp_j(\bs\te):=\diag_{k\in\{1,\cdots,q\}}\bigg\{\mathrm{e}^{\im (\l_j-\pi)\te_k/2}\left(\frac{j}{\mz p}\right)^{\!\!\te_k}\bigg\}\quad\mbox{and}\quad\mz p:=\exp\bigg(\frac{1}{m}\sum_{j=1}^m\log(j)\bigg),\]\normalsize
and, as $m$ tends to infinity, $\mz p\sim m/\mathrm{e}$. Observe that $\mcd(\bs d)$ is positive semidefinite since each summand of $\mcd$ is.  For $\kappa\in (0,1)$, define
\small\[\mdk(\bs d):=\sz\!\re\big[\mcp_j(\bs \te)\L_j(\do)^{-1}f_n(\l_j)\comp{\L_j(\do)^{-1}}\comp{\mcp_j(\bs \te)}\big]\]\normalsize
and
\small\[ \qz:=\sz\!\re\big[\mcp_j(\bs\te)G_0\comp{\mcp_j(\bs\te)}\big],\]\normalsize
where $[x]$ denotes the integer part of  $x$.
\small\begin{align}\label{auxl5}
\mdk(\bs d)&=\sz\!\re\big[\mcp_j(\bs \te)\L_j(\do)^{-1}\big(f(\l_j)+o_\P(n^{-\beta})\big)\comp{\L_j(\do)^{-1}}\comp{\mcp_j(\bs \te)}\big]\nonumber\\
&=\qz+\sz\re\big[\mcp_j(\bs \te)\L_j(\do)^{-1}o_\P(n^{-\beta})\comp{\L_j(\do)^{-1}}\comp{\mcp_j(\bs \te)}\big],
\end{align}\normalsize
where the last equality follows from lemma 5.4 in Shimotsu and Phillips (2005).
The $(r,s)$-th element of the third term on the RHS of \eqref{auxl5} is given by
\small\begin{align*}
\re\bigg[\sz  \left(\frac{j}{\mz p}\right)^{\!\!\te_r+\te_s}&\left(\frac{2\pi j}{n}\right)^{\!\!d^0_r+d^0_s}\mathrm{e}^{\im (\l_j-\pi)(\hat d_r-\hat d_s)/2}o_\P(n^{-\beta})\bigg]=\\
&=O(1)\left(\frac{m}{\mz p}\right)^{\!\!\te_r+\te_s}\left(\frac{m}{n}\right)^{\!\!d^0_r+d^0_s}o_\P(n^{-\beta})\,\sz \left(\frac{j}{m}\right)^{\!\!2(d^0_r+d^0_s)-(\widehat d_r+\widehat d_s)}\\
&=O(1)o_\P(1)\left(\frac{m}{\mz p}\right)^{\!\!\te_r+\te_s}O(1)=o_\P(1),
\end{align*}\normalsize
uniformly in $\bs \te\in\Theta_2$, since $\do\in\Omega_\beta$, $\beta\in(0,1)$ and Assumption \textbf{A}4, where the penultimate equality follows from lemma 5.4 in Shimotsu and Phillips (2005). Hence
\small\[\sup_{\Theta_2}\big\{\big|\det\{\mcd(\bs d)\}-\det\{\qz\}\big|\big\}=o_\P(1).\]\normalsize
The proof now follows viz a viz (with the obvious notational identification) from the proof of theorem 1 in Shimotsu (2007), p.294 (see the argument following equation (16)). We thus conclude that $P_2\rightarrow0$, as $n$ tends to infinity, and the proof is complete. \fim

\subsubsection*{Proof of Lemma \ref{lema2}:}
For fixed $r,s\in\{1,\cdots,q\}$, let $\mathscr{A}_{uv}:=\sum_{j=u}^v\mathscr A_j$ and $\mathscr{B}_{uv}:=\sum_{j=u}^v\mathscr B_j$, where
\small\begin{equation}\label{aux8}
\mathscr A_j:=\mathrm{e}^{\im(\l_j-\pi)(d_r^0-d_s^0)/2}\l_j^{d^0_r+d^0_s}\big[f_n^{rs}(\l_j) - \big(A(\l_j)\big)_{r\bs\cdot}I_{\bs\eps}(\l_j) \big(\comp{A(\l_j)}\big)_{\bs\cdot s}\big],
\end{equation}\normalsize
and
\small\begin{equation}\label{aux9}
\mathscr{B}_j:=\mathrm{e}^{\im(\l_j-\pi)(d_r^0-d_s^0)/2}\l_j^{d^0_r+d^0_s}\big(A(\l_j)\big)_{r\bs\cdot}I_{\bs\eps}(\l_j) \big(\comp{A(\l_j)}\big)_{\bs\cdot s}-G_0^{rs}.
\end{equation}\normalsize
Hence, for each $j$, $\mathscr A_j+\mathscr B_j=\mathrm{e}^{\im(\l_j-\pi)(d_r^0-d_s^0)/2}\l_j^{d^0_r+d^0_s}f_n^{rs}(\l_j)-G_0^{rs}$. For fixed $u\leq j\leq v$, we have
\small\[\E\big(|\mathscr A_j|\big)=\E\Big(\l_j^{d^0_r+d^0_s}\Big|f_n^{rs}(\l_j) - \big(A(\l_j)\big)_{r\bs\cdot}I_{\bs\eps}(\l_j) \big(\comp{A(\l_j)}\big)_{\bs\cdot s}\Big|\Big)=o(1),\]\normalsize
which yields $\displaystyle{\max_{r,s}}$ $\!\big\{\E\big(\big|\sum_{j=u}^v\mathscr{A}_j\big|\big)\big\} = o(v-u+1)$ and the result on $\mathscr A_{uv}$ follows.
As for $\mathscr B_j$, from the proof of lemma 1(a) in Shimotsu (2007) (notice that $\mathscr B_j$ does not depend on $f_n$) it follows that $\mathscr{B}_{uv}=o_\P(v)$ uniformly in $u$ and $v$ and the desired result follows.\fim

\subsubsection*{Proof of Theorem \ref{nonbeta}:}
From a careful inspection of the proof of Theorem \ref{consistency}, we observe that it suffices to show (with the same notation as in that proof)
\small\begin{equation}\label{reaux1}
\frac1m\sum_{j=1}^m\re\Big[M_j(\bs\te)\Lambda_j(\do)^{-1}f_n(\l_j)\comp{\Lambda_j(\do)^{-1}}\comp{M_j(\bs\te)}\Big] =\frac{1}{m}\sum_{j=1}^m\re\big[M_j(\bs\te) G_0 \comp{M_j(\bs\te)} \big]+o_\P(1),
\end{equation}\normalsize
uniformly in $\Theta_1$ and that $\mdk(\bs d)-\qz=o_\P(1)$, uniformly in $\Theta_2$. To show \eqref{reaux1}, notice that the $(r,s)$-th component of the LHS in \eqref{reaux1} is given by
\small\[\frac{1}{m}\sum_{j=1}^m\re\bigg[\mathrm{e}^{\im(\l_j-\pi)(\te_r-\te_s)/2}\left(\frac{j}{m}\right)^{\te_r+\te_s} \!\! f_n^{rs}(\l_j)\Big(\L_j^{(r)}(\do)\comp{\L_j^{(s)}(\do)}\Big)^{-1}\bigg].\]\normalsize
Summation by parts (see Zygmund, 2002, p.3) yields
\small\begin{eqnarray}\label{op1}
&&\hspace{-1cm}\sup_{\Theta_1}\bigg\{\bigg|\frac{1}{m}\sum_{j=1}^m\mathrm{e}^{\im(\l_j-\pi)(\te_r-\te_s)/2}\left(\frac{j}{m}\right)^{\te_r+\te_s} \!\! \Big[f_n^{rs}(\l_j)\Big(\L_j^{(r)}(\do)\comp{\L_j^{(s)}(\do)}\Big)^{\!\!-1}\!\!\!-G_0^{rs}\Big]\bigg|\bigg\}\leq\nonumber\\
&\leq&\frac{1}{m}\sum_{k=1}^{m-1}\sup_{\Theta_1}\bigg\{\bigg|\mathrm{e}^{\im(\l_k-\pi)(\te_r-\te_s)/2} \left(\frac{k}{m}\right)^{\te_r+\te_s}\hspace{-.3cm}-\mathrm{e}^{\im(\l_{k+1}-\pi)(\te_r-\te_s)/2} \left(\frac{k+1}{m}\right)^{\te_r+\te_s}\bigg|\bigg\}\times\nonumber\\
&&\times\,\,\,\bigg|\sum_{j=1}^k\Big[f_n^{rs}(\l_j)\Big(\L_j^{(r)}(\do)\comp{\L_j^{(s)}(\do)}\Big)^{\!\!-1}\!\!\!-G_0^{rs}\Big]\bigg| +\bigg|\frac{1}{m}\sum_{j=1}^m\Big[f_n^{rs}(\l_j)\Big(\L_j^{(r)}(\do)\comp{\L_j^{(s)}(\do)}\Big)^{\!\!-1}\!\!\!-G_0^{rs}\Big]\bigg|\nonumber\\
&\leq&\mz C\sum_{k=1}^{m-1}\left(\frac{k}{m}\right)^{\!\!2\epsilon}\frac{1}{k^2}\bigg|\sum_{j=1}^k\Big[f_n^{rs}(\l_j)\Big(\L_j^{(r)}(\do) \comp{\L_j^{(s)}(\do)}\Big)^{\!\!-1}\!\!\!-G_0^{rs}\Big]\bigg|+\nonumber\\
&&+\,\,\,\bigg|\frac{1}{m}\sum_{j=1}^m\Big[f_n^{rs}(\l_j)\Big(\L_j^{(r)}(\do)\comp{\L_j^{(s)}(\do)}\Big)^{\!\!-1}\!\!\!-G_0^{rs}\Big]\bigg|,
\end{eqnarray}\normalsize
\normalsize
where $0<\mz C<\infty$ is a constant. Now, from Lemma \ref{lema2},
\small\begin{align}\label{new1}
\sum_{k=1}^{m-1}\left(\frac{k}{m}\right)^{2\epsilon}\!\!\!\frac{1}{k^2}\bigg|\sum_{j=1}^k\Big[f_n^{rs}(\l_j)&\Big(\L_j^{(r)}(\do) \comp{\L_j^{(s)}(\do)}\Big)^{\!\!-1}\!\!\!-G_0^{rs}\Big]\bigg|\leq\sum_{k=1}^{m-1}\left(\frac{k}{m}\right)^{2\epsilon}\!\!\!\frac{1}{k^2}\Big(\big|\mathscr{A}_{1k}\big| +\big|\mathscr{B}_{1k}\big|\Big)\nonumber\\
&=\frac{1}{m^{2\epsilon}}\sum_{k=1}^{m-1}k^{2(\epsilon-1)}\big|\mathscr{A}_{1k}\big|+\frac{1}{m^{2\epsilon}}\sum_{k=1}^{m-1}k^{2(\epsilon-1)}o_\P(k).
\end{align}\normalsize
The RHS of \eqref{new1} is $o_\P(1)$ uniformly in $(r,s)$ by lemma 1 in Robinson (1995b), Lemma \ref{lema2} and Chebyshev's inequality since $\E\big(m^{-2\epsilon}\sum_{k=1}^{m-1}k^{2(\epsilon-1)}\big|\mathscr{A}_{1k}\big|\big)=o(1)$, uniformly in $(r,s)$. The other term in \eqref{op1} is also $o_\P(1)$, uniformly in $(r,s)$, by the same argument and, hence, \eqref{reaux1} follows. On the other hand, the $(r,s)$-th element of $\mdk(\bs d)-\qz$ is given by
\small
\begin{align}
&\sz\!\re\bigg[\mathrm{e}^{\im(\l_j-\pi)(\te_r-\te_s)/2}\left(\frac{j}{\mz p}\right)^{\te_r+\te_s} \!\! \Big[f_n^{rs}(\l_j)\Big(\L_j^{(r)}(\do)\comp{\L_j^{(s)}(\do)}\Big)^{\!\!-1}\!\!\!-G_0^{rs}\Big]\bigg]=\nonumber\\
&\ =\left(\frac{m}{\mz p}\right)^{\te_r+\te_s} \!\!\re\bigg[\sz\mathrm{e}^{\im(\l_j-\pi)(\te_r-\te_s)/2}\left(\frac{j}{m}\right)^{\te_r+\te_s} \!\!\Big[f_n^{rs}(\l_j)\Big(\L_j^{(r)}(\do)\comp{\L_j^{(s)}(\do)}\Big)^{\!\!-1}\!\!\!-G_0^{rs}\Big]\bigg]=o_\P(1),\nonumber
\end{align}
\normalsize
uniformly in $\bs\te\in\Theta_2$, where the last equality is derived similarly to \eqref{op1} from summation by parts and lemma 5.4 in Shimotsu and Phillips (2005). This completes the proof.\fim

\subsubsection*{Proof of Lemma \ref{lemma1b2}:}
\rm{(a)} For $r,s\in\{1,\cdots,q\}$ fixed, we have
\small\begin{align*}
\sum_{j=1}^v\mathrm{e}^{\im(\l_j-\pi)(d_r^0-d_s^0)/2}&\l_j^{d_r^0+d_s^0}\Big[f_n^{rs}(\l_j)- \big(A(\l_j)\big)_{r\bs\cdot}I_{\bs\eps}(\l_j)\big(\comp{A(\l_j)}\big)_{\bs\cdot s}\Big]\leq\\
&\leq\Big[1+O\big((\l_j-\pi)(d_r^0-d_s^0)/2\big)\Big]\sum_{j=1}^v\l_j^{d_r^0+d_s^0}\Big[f_n^{rs}(\l_j)- \big(A(\l_j)\big)_{r\bs\cdot}I_{\bs\eps}(\l_j)\big(\comp{A(\l_j)}\big)_{\bs\cdot s}\Big]\\
&\leq O(1)\max\{1,n^{|d_r^0+d_s^0|}\}\max_{1\leq v\leq m}\bigg\{\sum_{j=1}^v\Big[f_n^{rs}(\l_j)- \big(A(\l_j)\big)_{r\bs\cdot}I_{\bs\eps}(\l_j)\big(\comp{A(\l_j)}\big)_{\bs\cdot s}\Big]\bigg\}\\
&=o_\P\bigg(\frac{m^2}n\bigg)=o_\P\bigg(\frac{\sqrt m}{\log(m)}\bigg),
\end{align*}\normalsize
where the last equality follows from \textbf B4, and \eqref{lema1b2a} follows.

\noindent\rm{(b)}  Rewrite the argument of the summation in \eqref{lema1b2b} as $\mathscr{A}_j+\mathscr{B}_j+\mathscr{C}_j$, where
\small\begin{align*}
\mathscr{A}_j&:=\mathrm{e}^{\im(\l_j-\pi)(d_r^0-d_s^0)/2}\l_j^{d_r^0+d_s^0}\Big[f_n^{rs}(\l_j)- \big(A(\l_j)\big)_{r\bs\cdot}I_{\bs\eps}(\l_j)\big(\comp{A(\l_j)}\big)_{\bs\cdot s}\Big],\\
\mathscr{B}_j&:=\mathrm{e}^{\im(\l_j-\pi)(d_r^0-d_s^0)/2}\l_j^{d_r^0+d_s^0}\Big[ \big(A(\l_j)\big)_{r\bs\cdot}I_{\bs\eps}(\l_j)\big(\comp{A(\l_j)}\big)_{\bs\cdot s}-f_{rs}(\l_j)\Big],\\
\mathscr{C}_j&:=\mathrm{e}^{\im(\l_j-\pi)(d_r^0-d_s^0)/2}\l_j^{d_r^0+d_s^0}f_{rs}(\l_j)-G_0^{rs}.
\end{align*}\normalsize
Part \rm{(a)} yields $\max_{r,s}\Big\{\!\!\sum_{j=1}^v\big|\mathscr{A}_j\big|\Big\} =o_\P\Big(m^{1/2}\big(\log(m)\big)^{-1}\Big)$, while, from the proof of lemma 1(b2) in Shimotsu (2007), we obtain $\max_{r,s}\Big\{\!\!\sum_{j=1}^v\big|\mathscr{B}_j\big|\Big\} =O_\P\big(m^{1/2}\log(m)\big)$. Assumption \textbf{B}1 implies $\max_{r,s}\Big\{\!\!\sum_{j=1}^v\big|\mathscr{C}_j\big|\Big\} =O\big(m^{\alpha+1}n^{-\alpha}\big)$. The result now follows by noticing that $m^{1/2}\log(m)^{-1}=O\big(m^{1/2}\log(m)\big)$.\fim

\subsubsection*{Proof of Theorem \ref{norm}:}
The idea of the proof is similar to that of Lobato (1999) with similar adaptations as in Shimotsu (2007). By hypothesis,
\small\[\bs0=\frac{\partial S(\bs d)}{\partial\bs d}\bigg|_{\d}=\frac{\partial S(\bs d)}{\partial\bs d}\bigg|_{\do}+\bigg(\frac{\partial^2S(\bs d)}{\partial \bs d\partial\bs{d}^{\,\prime}}\bigg|_{\od}\bigg)(\d-\do),\]\normalsize
with probability tending to 1, as $n$ tends to infinity, for some $\od$ such that $\|\od-\do\|_{\infty}\leq\|\d-\do\|_{\infty}$. We observe that $\d$ has the stated limiting distribution if
\small\begin{equation}\label{M1}
\sqrt{m}\,\frac{\partial S(\bs d)}{\partial\bs d}\bigg|_{\do}\overset{d}{-\!\!\!\longrightarrow} N(0,\Sigma)
\end{equation}\normalsize
and
\small\begin{equation}\label{M2}
\frac{\partial^2S(\bs d)}{\partial \bs d\partial\bs{d}^{\,\prime}}\bigg|_{\od}\overset{\P}{-\!\!\!\longrightarrow}\Sigma.
\end{equation}\normalsize
We shall prove \eqref{M1} first. Observe that, for $r\in\{1,\cdots,q\}$,
\small\[\sqrt{m}\,\frac{\partial S(\bs d)}{\partial d_r}=-\frac{2}{\sqrt{m}}\sum_{j=1}^m \log(\l_j)+\tr\bigg[\G{\bs d}^{-1}\sqrt{m}\,\frac{\partial\G{\bs d}}{\partial d_r}\bigg].\]\normalsize
Let $\ir $ denote a $q\times q$ matrix whose $(r,r)$-th element is 1 and all other elements are zero.
Define a function $\varphi:(0,\infty)\rightarrow\C$ by
\small\begin{equation}\label{l(x)}
\varphi(x):=\log(x)+\im\bigg(\frac{x-\pi}{2}\bigg).
\end{equation}\normalsize
Since $\L_j(\bs d)^{-1}=\displaystyle{\diag_{k\in\{1,\cdots,q\}}}\big\{\l_j^{d_k}\mathrm{e}^{\im(\l_j-\pi)d_k/2}\big\}$ and $\re\big[(a+\im b)(c+\im d) \big] =ac-bd$, we can write
\small
\begin{align}\label{n1}
\sqrt{m}\,\frac{\partial\G{\bs d}}{\partial d_r}\bigg|_{\do}&=\frac{1}{\sqrt{m}}\sum_{j=1}^m\re\Big[\varphi(\l_j)\L_j(\do)^{-1}\ir f_n(\l_j)\comp{\L_j(\do)^{-1}} \Big]+\nonumber\\
&\qquad+\frac{1}{\sqrt{m}}\sum_{j=1}^m\re\Big[\overline{\varphi(\l_j)}\L_j(\do)^{-1}f_n(\l_j)\ir \comp{\L_j(\do)^{-1}} \Big]\nonumber\\
&=\frac{1}{\sqrt{m}}\sum_{j=1}^m\log(\l_j)\re\Big[\L_j(\do)^{-1}\big(\ir f_n(\l_j)+f_n(\l_j)\ir \big) \comp{\L_j(\do)^{-1}}\Big] +\nonumber\\
&\qquad+\frac{1}{\sqrt{m}}\sum_{j=1}^m\bigg[\frac{\l_j-\pi}{2}\bigg]\pim\Big[\L_j(\do)^{-1}\big(-\ir f_n(\l_j)+ f_n(\l_j)\ir \big)\comp{\L_j(\do)^{-1}}\Big],\nonumber\\
&:=\mathscr{H}_1(r)+\mathscr{H}_2(r).
\end{align}
\normalsize
Therefore,  for $\eta$ an arbitrary vector in $\R^q$, from \eqref{n1} we obtain
\small\begin{align*}
\eta^\prime\sqrt{m}\,&\frac{\partial S(\bs d)}{\partial\bs d}\bigg|_{\do}=\sk\eta_k\sqrt{m}\,\frac{\partial S(\bs d)}{\partial d_k}\bigg|_{\do}=\\
&=\sk\eta_k\bigg[-\frac{2}{\sqrt{m}}\sum_{j=1}^m\log(\l_j)+\tr\big[\G{\do}^{-1}\mathscr{H}_1(k)\big]\bigg] + \sk\eta_k\tr\big[\G{\do}^{-1}\mathscr{H}_2(k)\big],\\
& :=\mathscr{R}_1+\mathscr{R}_2.
\end{align*}
\normalsize
We analyze $\mathscr R_1$ first. By letting
\small\[a_j:=\log(\l_j)-\frac{1}{m}\sum_{k=1}^m\log(\l_k)=\log(j)-\frac{1}{m}\sum_{k=1}^m\log(k)=O\big(\log(m)\big),\]\normalsize
we can write
\small
\begin{align}\label{n2}
-\frac{2}{\sqrt{m}}\sum_{j=1}^m\log(\l_j)+\tr\big[\G{\do}^{-1}&\mathscr{H}_1(k)\big]=\tr\bigg[\G{\do}^{-1}\bigg(\mathscr H_1(k)-\frac{2}{\sqrt{m}}\sum_{j=1}^m\log(\l_j)\G{\do}\mathrm{I}_{(k)}\bigg)\bigg]\nonumber\\
&=\tr\bigg[\G{\do}^{-1}\frac{2}{\sqrt{m}}\sum_{j=1}^ma_j\re\Big[\lfl\Big]\mathrm{I}_{(k)}\bigg].
\end{align}
\normalsize
By Lemma \ref{lemma1b2}\rm{(b)}, \eqref{n2} can be written as
\small\begin{equation}
\Big[(G_0^{-1})_{k\bs\cdot}+o_\P(1)\Big]\frac{2}{\sqrt{m}}\sum_{j=1}^ma_j\Big(\re\Big[\lfl\Big]\Big)_{\bs\cdot k}.\nonumber
\end{equation}\normalsize
Now, by Lemma \ref{lemma1b2}(a),
\small
\begin{align*}
\bigg(\sum_{j=1}^m&a_j\L_j(\do)^{-1}\Big(f_n(\l_j)-A(\l_j)I_{\bs\eps}(\l_j)\comp{A(\l_j)}\Big)\comp{\L_j(\do)^{-1}}\bigg)_{rs}\leq\\
&\leq O\big(\log(m)\big)\!\max_{v=1,\cdots,m}\!\!\bigg\{\sum_{j=1}^v\mathrm{e}^{\im(\l_j-\pi)(d_r^0-d_s^0)/2}\l_j^{d_r^0+d_s^0} \Big(f_n^{rs}(\l_j)-\big(A(\l_j)\big)_{r\bs\cdot}I_{\bs\eps}(\l_j)\big(\comp{A(\l_j)}\big)_{\bs\cdot s}\Big)\bigg\}\\
&=O\big(\log(m)\big)o_\P\bigg(\frac{\sqrt{m}}{\log(m)}\bigg)\,\,=\,\,o_\P(\sqrt{m}),
\end{align*}
\normalsize
uniformly in $r,s\in\{1,\cdots, q\}$. Therefore,
\small\begin{align}\label{n4}
\frac{1}{\sqrt{m}}\sum_{j=1}^m&a_j\lfl=\nonumber\\
&=\frac{1}{\sqrt{m}}\sum_{j=1}^ma_j \L_j(\do)^{-1}A(\l_j)I_{\bs\eps}(\l_j)\comp{A(\l_j)}\comp{\L_j(\do)^{-1}} + \frac{1}{\sqrt{m}}\,o_\P(\sqrt{m})\nonumber\\
&=\frac{1}{\sqrt{m}}\sum_{j=1}^ma_j \Big[\L_j(\do)^{-1}A(\l_j)I_{\bs\eps}(\l_j)\comp{A(\l_j)}\comp{\L_j(\do)^{-1}}-G_0\Big]+o_\P(1),
\end{align}\normalsize
where the last equality follows from $\sum_{j=1}^m a_j=0$. The proof of \eqref{M1} now follows viz a viz from the proof of theorem 2 in Shimotsu (2007) p.296, by  noticing that, with the appropriate notational identification, \eqref{n4} is (21) in the aforementioned theorem.

We move to prove \eqref{M2}. For fixed $\delta>0$, let $\bs\te:=\bs d-\do$ and define
\small\[\mathcal M:=\Big\{\bs d:\log(n)^4\|\bs d -\do\|_{\infty}<\delta\Big\}=\Big\{\bs \te :\log(n)^4\|\bs\te\|_{\infty}<\delta\Big\}.\]\normalsize
First we show that $\P\big(\od\in \mathcal M\big)\rightarrow 1$, as $n\rightarrow\infty$. Assuming the same notation  as in the proof of Theorem \ref{consistency},  recall the decomposition of $L(\bs d)=S(\bs d)-S(\do)$ given in expression \eqref{aux1}. By applying the same argument as in the proof of Theorem \ref{consistency}, we first obtain
\small\[\inf_{\Theta_1\setminus \mathcal M}\big\{\mz R(\bs d)\big\}\geq\delta^2\log(n)^8,\]\normalsize
and upon applying  Lemma \ref{lemma1b2}, we obtain
\small\[\sup_{\Theta_1}\Big\{\big|\mz A(\bs d)-h(\bs d)\big|\Big\}=O_\P\bigg(\frac{m^\alpha}{n^\alpha}+\frac{\log(m)}{m^{2\epsilon}}+\frac{m}{n}\bigg).\]\normalsize
It follows, uniformly in $\Theta_1$ (cf. Shimotsu, 2007, p.300), that
\small\begin{align*}
\log\big(\mz A(\bs d)\big)&-\log\big(\mz B(\bs d)\big)\geq \log\big(h(\bs d)+o_\P\big(\log(n)^{-8}\big)\big)-\log\big(h(\bs d)\big)=o_\P\big(\log(n)^{-8}\big)\\
\log\big(\mz A(\do)\big)&-\log\big(\mz B(\do)\big)=\log\big(h(\do)+o_\P\big(\log(n)^{-8}\big)\big)-\log\big(h(\do)\big)=o_\P\big(\log(n)^{-8}\big).
\end{align*}\normalsize
Hence, $\P\big(\inf_{\Theta_1\setminus\mathcal M}L(\bs d)\leq 0\big)\rightarrow0$ and $\P\big(\od\in \mathcal M\big)\rightarrow 1$, as $n\rightarrow\infty$, follows.

Now, observe that
\small\[\frac{\partial^2 S(\bs d)}{\partial d_r\partial d_s}=\tr\bigg[-\G{\bs d}^{-1}\,\frac{\partial\G{\bs d}}{\partial d_r}\,\G{\bs d}^{-1}\frac{\partial\G{\bs d}}{\partial d_s}+\G{\bs d}^{-1}\frac{\partial^2 \G{\bs d}}{\partial d_r\partial d_s}\bigg].\]\normalsize
The derivatives of $\G{\bs d}$ are given by
\small
\begin{align}\label{h1}
\frac{\partial\G{\bs d}}{\partial d_r}=\frac{1}{m}\sum_{j=1}^m&\,\re\Big[\varphi(\l_j)\ir \L_j(\bs d)^{-1}f_n(\l_j)\comp{\L_j(\bs d)^{-1}}\Big]+\frac{1}{m}\sum_{j=1}^m\re\Big[\overline{\varphi(\l_j)}\L_j(\bs d)^{-1}f_n(\l_j)\comp{\L_j(\bs d)^{-1}}\ir \Big],
\end{align}
\normalsize
and
\small\begin{align*}
\frac{\partial^2 \G{\bs d}}{\partial d_r\partial d_s}=&\,\frac{1}{m}\sum_{j=1}^m\re\Big[\varphi(\l_j)^2\ir\is \L_j(\bs d)^{-1}f_n(\l_j)\comp{\L_j(\bs d)^{-1}}\Big]+\\
&+\frac{1}{m}\sum_{j=1}^m\re\Big[\big|\overline{\varphi(\l_j)}\big|^2\ir \L_j(\bs d)^{-1}f_n(\l_j)\comp{\L_j(\bs d)^{-1}} \is\Big]+\\
&+\frac{1}{m}\sum_{j=1}^m\re\Big[\big|\overline{\varphi(\l_j)}\big|^2\is\L_j(\bs d)^{-1}f_n(\l_j)\comp{\L_j(\bs d)^{-1}} \ir \Big]+\\
&+\frac{1}{m}\sum_{j=1}^m\re\Big[\overline{\varphi(\l_j)}^{\,2}\L_j(\bs d)^{-1}f_n(\l_j)\comp{\L_j(\bs d)^{-1}} \ir\is \Big],
\end{align*}\normalsize
where $\varphi$ is given by \eqref{l(x)}. For $k=0,1,2$, let
\small\begin{align*}
\mathcal R_k(\bs d)&=\frac{1}{m}\sum_{j=1}^m\log(\l_j)^k\re\Big[\L_j(\bs d)^{-1}f_n(\l_j)\comp{\L_j(\bs d)^{-1}} \Big]\\
\mathcal I_k(\bs d)&=\frac{1}{m}\sum_{j=1}^m\log(\l_j)^k\pim\Big[\L_j(\bs d)^{-1}f_n(\l_j)\comp{\L_j(\bs d)^{-1}}\Big],
\end{align*}\normalsize
so that, we can write
\small
\begin{equation}\label{deriv1}
\frac{\partial\G{\bs d}}{\partial d_r}=\ir \mathcal R_1(\bs d)+\mathcal R_1(\bs d)\ir +\frac{\pi}{2}\Big(\ir \mathcal I_0(\bs d)-\mathcal I_0(\bs d)\ir \Big)+o_\P\left(\frac{1}{\log(n)}\right),
\end{equation}
\normalsize
and
\small
\begin{align}\label{deriv2}
\frac{\partial^2\G{\bs d}}{\partial d_r\partial d_s}&=\frac{\pi^2}{4}\Big[\ir\is \mathcal R_0(\bs d) +\ir \mathcal R_0(\bs d)\is+\is\mathcal R_0(\bs d)\ir  +\mathcal R_0(\bs d)\is\ir \Big]+\nonumber\\
&\hspace{1cm}+\pi \Big[\ir\is \mathcal I_1(\bs d)+\mathcal I_1(\bs d)\ir\is \Big] +\ir\is \mathcal R_2(\bs d) + \ir \mathcal R_2(\bs d)\is+\nonumber\\
&\hspace{2cm}+\is\mathcal R_2(\bs d)\ir  +\mathcal R_2(\bs d)\is\ir + o_\P(1).
\end{align}
\normalsize
The order of the remainder term in \eqref{deriv1} is obtained as follows. Rewrite the first term on the RHS of \eqref{h1} as
\small
\begin{align*}
\frac{1}{m}&\sum_{j=1}^m\re\Big[\varphi(\l_j)\ir \Big(\re\Big[\L_j(\bs d)^{-1}f_n(\l_j)\comp{\L_j(\bs d)^{-1}}\Big]+\im\,\pim\Big[\L_j(\bs d)^{-1}f_n(\l_j)\comp{\L_j(\bs d)^{-1}}\Big]\Big) \Big]=\\
&=\frac{1}{m}\sum_{j=1}^m\re\bigg[\log(\l_j)\ir \re\Big[\L_j(\bs d)^{-1}f_n(\l_j)\comp{\L_j(\bs d)^{-1}}\Big]+ \im\log(\l_j)\ir \pim\Big[\L_j(\bs d)^{-1}f_n(\l_j)\comp{\L_j(\bs d)^{-1}}\Big]+\\
&\hspace{.5cm}+\im\bigg(\frac{\l_j-\pi}{2}\bigg)\ir \re\Big[\L_j(\bs d)^{-1}f_n(\l_j)\comp{\L_j(\bs d)^{-1}}\Big] -\bigg(\frac{\l_j-\pi}{2}\bigg)\ir \pim\Big[\L_j(\bs d)^{-1}f_n(\l_j)\comp{\L_j(\bs d)^{-1}}\Big]\bigg]\\
&=\frac{1}{m}\sum_{j=1}^m\bigg[\log(\l_j)\ir \re\Big[\L_j(\bs d)^{-1}f_n(\l_j)\comp{\L_j(\bs d)^{-1}}\Big]-\bigg(\frac{\l_j-\pi}{2}\bigg)\ir \pim\Big[\L_j(\bs d)^{-1}f_n(\l_j)\comp{\L_j(\bs d)^{-1}}\Big]\bigg].
\end{align*}
\normalsize
By summation by parts, Lemma \ref{lemma1b2} and assumption \textbf B4
\small\begin{align*}
\frac{1}{m}\sum_{j=1}^m\l_j\Big[\L_j(\bs d)^{-1}f_n(\l_j)\comp{\L_j(\bs d)^{-1}}\Big]&\leq \frac{1}{m}\sum_{j=1}^{m-1}\big|\l_j-\l_{j+1}\big|\bigg\| \sum_{k=1}^j\L_j(\bs d)^{-1}f_n(\l_j)\comp{\L_j(\bs d)^{-1}}\bigg\|_{\infty}+\\
&\hspace{1cm}+\frac{\l_m}{m}\bigg\| \sum_{k=1}^j\L_j(\bs d)^{-1}f_n(\l_j)\comp{\L_j(\bs d)^{-1}}\bigg\|_{\infty}\\
&\hspace{-3cm}\leq\frac{1}{m}\frac{m-1}{n}\bigg[O_\P\left((m-1)^{{1/2}}\log(m-1)+\frac{(m-1)^{\alpha+1}}{n^\alpha}\right) +O\left(\frac{1}{m}\right)\bigg]+\\
&\hspace{-2cm}+O\left(\frac{1}{n}\right)\bigg[O_\P\left(m^{{1/2}}\log(m)+\frac{m^{\alpha+1}}{n^\alpha}\right) +O\left(\frac{1}{m}\right)\bigg]\\
&\hspace{-3cm}=O_\P\left(\frac{m^{{1/2}}\log(m)}{n}+\frac{m^{\alpha}}{n^{\alpha}}\right)+ O\left(\frac{1}{mn}\right)\,\,=\,\,o_\P\left(\frac{1}{\log(n)}\right),
\end{align*}\normalsize
where the last equality follows from
\small\[\frac{m^{1/2}\log(m)}{n}\,\log(n)=\frac{m^{1/2}}{n^{1/2}}\,\frac{\log(m)}{n^{1/4}}\frac{\log(n)}{n^{1/4}}\longrightarrow 0,\]\normalsize
as $n$ tends to infinity, by assumption \textbf B4. The other term in \eqref{h1} is dealt analogously. The remainder term in \eqref{deriv2} involves
\small\[\frac{1}{m}\sum_{j=1}^m\rho(\l_j)\L_j(\bs d)^{-1}f_n(\l_j)\comp{\L_j(\bs d)^{-1}},\]\normalsize
for $\rho(\l_j)$ proportional to $\l_j,\l_j^2$ and $\l_j\log(\l_j)$. The order of the terms proportional to $\l_j$ has already been obtained, while the one proportional to $\l_j^2$ is dealt analogously since $\l_j^2=O(\l_j)$. The term proportional to $\l_j\log(\l_j)$ is $o_\P(1)$ since, by summation by parts, Lemma \ref{lemma1b2} and assumption \textbf B4,
\small
\begin{align*}
\bigg\|\frac{1}{m}\sum_{j=1}^m\l_j\log(\l_j)&\Big[\L_j(\bs d)^{-1}f_n(\l_j)\comp{\L_j(\bs d)^{-1}}\Big]\bigg\|_\infty\leq\\
&\leq\frac{1}{m}\sum_{j=1}^{m-1}\Big|\l_j\log(\l_j)-\l_{j+1}\log(\l_{j+1})\Big|\bigg\|\sum_{k=1}^j\L_k(\bs d)^{-1}f_n(\l_j)\comp{\L_k(\bs d)^{-1}}\bigg\|_\infty+\\
&\hspace{1cm}+\frac{\l_m\log(\l_m)}{m}\bigg\|\sum_{j=1}^m\L_j(\bs d)^{-1}f_n(\l_j)\comp{\L_j(\bs d)^{-1}}\bigg\|_\infty\\
&=\frac{1}{m}\,O_\P\left((m-1)^{1/2}\log(m-1)+\frac{(m-1)^{\alpha+1}}{n^\alpha}\right)o(1)+o_\P(1)\ =\ o_\P(1).
\end{align*}
\normalsize
In order to finish the proof, it suffices to show that,
\small\begin{equation}\label{h2}
\mathcal R_k(\bs d)=G_0\,\frac{1}{m}\sum_{j=1}^m\log(\l_j)^k+o_\P\big(\log(n)^{k-2}\big)
\end{equation}\normalsize
and
\small\begin{equation}\label{h3}
\mathcal I_k(\bs d)=o_\P\big(\log(n)^{k-2}\big),
\end{equation}\normalsize
uniformly in $\bs d\in \mathcal M$. Indeed, if \eqref{h2} and \eqref{h3} hold, upon defining for a matrix $M$,
\small\[T_1(M,r):=\ir M+M\ir,\qquad T_2(M,r,s):=\ir\is M+\ir M\is+\is M\ir+M\ir\is\]\normalsize
and
\small\[T_3(M,r,s):=-\ir\is M+\ir M\is+\is M\ir-M\ir\is,\]\normalsize
it follows that (cf. Shimotsu, 2007, p.301)
\small
\begin{align*}
&\G{\od}^{-1}=G_0^{-1}+o_\P\big(\log(n)^{-2}\big),\qquad\frac{\partial \G{\od}}{\partial d_r}=\frac{1}{m}\sum_{j=1}^m\log(\l_j)T_1(G_0,r)+o_\P\big(\log(n)^{-1}\big),\\
&\qquad\text{ and }\quad\frac{\partial^2 \G{\od}}{\partial d_r\partial d_s}=\frac{1}{m}\sum_{j=1}^m\log(\l_j)^2T_2(G_0,r,s)+\frac{\pi^2}{4}\,T_3(G_0,r,s)+o_\P(1).
\end{align*}
\normalsize
Since $\tr\Big[G_0^{-1}T_1(G_0,r)G_0^{-1}T_1(G_0,s)\Big]=\tr\big[G_0^{-1}T_2(G_0,r,s)\big]$ and
\small\[\frac{1}{m}\sum_{j=1}^m\log(\l_j)^2-\bigg(\frac{1}{m}\sum_{j=1}^m\log(\l_j)\bigg)^2\longrightarrow 1,\]\normalsize
we obtain
\small\[\frac{\partial^2 S(\bs d)}{\partial d_r\partial d_s}=\tr \left[G_0^{-1}T_2(G_0,r,s)+\frac{\pi^2}{4}\,G_0^{-1}T_3(G_0,r,s)\right]+o_\P(1),\]\normalsize
from where \eqref{M2} follows. We proceed to show \eqref{h2} and \eqref{h3}. For $k=0,1,2$, let
\small\[\mz F_k(\bs \te):=\frac{1}{m}\sum_{j=1}^m\log(\l_j)^k\L_j(\bs \te)^{-1}G_0\comp{\L_j(\bs \te)^{-1}}.\]\normalsize
Then, \eqref{h2} and \eqref{h3} follow if
\small\begin{equation}\label{h5}
\sup_{\bs d\in \mathcal M}\bigg\{\bigg\|\frac{1}{m}\sum_{j=1}^m\log(\l_j)^k\L_j(\bs d)^{-1}f_n(\l_j)\comp{\L_j(\bs d)^{-1}}-\mz F_k(\bs \te)\bigg\|_\infty\bigg\} =o_\P\big(\log(n)^{k-2}\big),
\end{equation}\normalsize
and
\small\begin{equation}\label{h6}
\sup_{\bs d\in \mathcal M}\bigg\{\bigg\|\mz F_k(\bs \te)-G_0\frac{1}{m}\sum_{j=1}^m\log(\l_j)^k\bigg\|_\infty\bigg\} =o\big(\log(n)^{k-2}\big).
\end{equation}\normalsize
Following Shimotsu (2007), p.302, notice that, by applying \eqref{4.15A}, we can rewrite \eqref{h5} as
\small\[\sup_{\bs d\in \mathcal M}\bigg\{\bigg\|\frac{1}{m}\sum_{j=1}^m\log(\l_j)^k\L_j(\bs \te)^{-1}\Big[\lfl-G_0\Big]\comp{\L_j(\bs \te)^{-1}}\bigg\|_\infty\bigg\}.\]\normalsize
Define $b_j(\bs\te;k):=\log(\l_j)^k\mathrm{e}^{\im(\l_j-\pi)(\te_r-\te_s)/2}\l_j^{\te_r+\te_s}$, for $k=0,1,2$. Then, by omitting the supremum, the $(r,s)$-th element of \eqref{h5} is equal to
\small\begin{align}\label{h7}
\bigg|\frac{1}{m}\sum_{j=1}^mb_j(\bs\te;k)&\Big[\mathrm{e}^{\im(\l_j-\pi)(d_r^0-d_s^0)/2}\l_j^{d^0_r+d^0_s}f_n^{rs}(\l_j)-G_0^{rs}\Big] \bigg|\leq\nonumber\\
&\leq\frac{1}{m}\sum_{j=1}^{m-1}\Big|b_j(\bs\te;k)-b_{j+1}(\bs\te;k) \Big| \bigg|\sum_{l=1}^{j}\mathrm{e}^{\im(\l_l-\pi)(d_r^0-d_s^0)/2}\l_l^{d^0_r+d^0_s} f_n^{rs}(\l_l)-G_0^{rs}\bigg|\nonumber\\
&\hspace{1cm}+\frac{b_m(\bs \te;k)}{m}\,\bigg|\sum_{j=1}^m\mathrm{e}^{\im(\l_j-\pi)(d_r^0-d_s^0)/2}\l_j^{d^0_r+d^0_s} f_n^{rs}(\l_j)-G_0^{rs}\bigg|,
\end{align}\normalsize
where the inequality follows from summation by parts. Now, since
\small\[b_j(\bs\te;k)-b_{j+1}(\bs\te;k) =O\left(\frac{\log(n)^k}{j}\right)\quad\mbox{ and }\quad b_m(\bs \te;k)=O\big(\log(n)^{k}\big),\]\normalsize
uniformly in $\bs \te\in\mathcal M$, for any $k=0,1,2$, it follows by Lemma \ref{lemma1b2} and Remark \ref{rmk} that \eqref{h7} can be rewritten as
\small\begin{align*}
O\left(\frac{\log(n)^k}{m}\right)&\frac{1}{m}\,O_\P\left(\frac{m^{\alpha+1}}{n^\alpha}+m^{1/2}\log(m)\right)+\\
&\hspace{1cm}+ O\big(\log(n)^k\big)\frac{1}{m}\,O_\P\left(\frac{m^{\alpha+1}}{n^\alpha}+m^{1/2}\log(m)\right) =o_\P\big(\log(n)^{k-2}\big),
\end{align*}\normalsize
where the last equality follows from assumption \textbf B4 (see also Remark \ref{rmk}), because
\small\begin{align*}
\log(n)^2\frac{1}{m}\,O_\P\left(\frac{m^{\alpha+1}}{n^\alpha}+m^{1/2}\log(m)\right)&=
\bigg[\frac{\log(n)^2}{m^{1/2}\log(m)}+\frac{\log(m)}{m^{1/4}}\,\frac{\log(n)^2}{m^{1/4}}\bigg]O_\P(1)=o_\P(1).
\end{align*}\normalsize
The other term is dealt analogously, so that \eqref{h5} follows. As for \eqref{h6}, the result follows from the proof of theorem 2, p.302, in Shimotsu (2007)  (notice that it does not depend on $f_n$). This completes the proof.\fim
\subsubsection*{Proof of Corollary \ref{tapered.cons}}
The proof follows the same lines as the proof of lemma 1(a) in Shimotsu (2007) p.308 in view of $I_T(\lambda;n)=O\big(I_n(\lambda)\big)$. \fim

\subsubsection*{Proof of Corollary \ref{tapered.an}}

 By carefully inspecting  the proof of Theorem \ref{norm}, we notice that suffices to show that part (a) and (b) of Lemma \ref{lema1b2b} hold for the result of Theorem \ref{norm} to hold. Part (a) and (b) of Lemma \ref{lema1b2b} are proven following the same lines as the proof of lemma 1(b1) and lemma 1(b2) in Shimotsu (2007), in view of $I_T(\lambda;n)=O\big(I_n(\lambda)\big)$. \fim


\small
\begin{thebibliography}{20}
\bibitem{c}  Chiriac, R. and Voev, V. (2011). ``Modeling and Forecasting Multivariate Realized Volatility''. {\slshape Journal of Applied Econometrics}, \textbf{26}(6), 922-947.
\bibitem{dah2} Dahlhaus, R. (1983). ``Spectral Analysis with Tapered Data''. {\slshape Journal of Time Series Analysis}, \textbf 4, 163-175.
\bibitem{ft} Fox, R. and Taqqu, M.S. (1986). ``Large-Sample Properties of Parameter Estimates for Strongly Dependent Stationary Gaussian Time Series''. {\slshape The Annals of Statistics}, \textbf{14}(2), 517-532.
\bibitem{fr} Fryzlewicz, P.,  Nason, G.P. and von Sachs, R. (2008). ``A Wavelet-Fisz Approach to Spectrum Estimation''. {\slshape Journal of Time Series Analysis}, \textbf{29}(5), 868-880.
\bibitem{gs} Giraitis, L. and Surgailis, D. (1990).``A Central Limit Theorem for Quadratic Forms in Strongly Dependent Linear Variables and its Application to Asymptotic Normality of Whittle's Estimate''. {\slshape Probability Theory and Related Fields}, \textbf{86}, 87-104.
\bibitem{g} Grenander, U. (1951). ``On Empirical Spectral Analysis of Stochastic Process''. {\slshape Arkiv f\"or Matematik}, \textbf{1}, 197-277.
\bibitem{hur2} Hurvich, C.M. and Beltr\~ao, K.I. (1993). ``Asymptotics for the Low-Frequency Ordinates of the Periodogram of a Long-Memory Time Series''. {\slshape Journal of Time Series Analysis}, \textbf{14}(5), 455-472.
\bibitem{hur3} Hurvich, C.M. and B.K. Ray (1995). ``Estimation of the Memory Parameter for Nonstationary or Noninvertible Fractionally Integrated Processes''. {\slshape Journal of Time Series Analysis}, \textbf{16}(1), 17-42.
\bibitem{k} K\"unsch, H. (1987). ``Statistical Aspects of Self-Similar Processes''. In Prokhorov, Yu. and Sazanov, V.V. (Eds.), {\slshape Proceedings of the First World Congress of the Bernoulli Society}. Utrecht: VNU Science Press, 67-74.
\bibitem{olv} Olbermann, B.P., Lopes, S.R.C. and Reisen, V.A. (2006). ``Invariance of the First Difference in ARFIMA Models''. {\slshape Computational Statistics}, \textbf{21}(3), 445-461.
\bibitem{lo}  Lobato, I.N. (1999). ``A Semiparametric Two-Step Estimator in a Multivariate Long Memory Model''. {\slshape Journal of Econometrics}, \textbf{90}, 129-153.
\bibitem{lu}  Luce\~no, A. (1996). ``A Fast Likelihood Approximation for Vector General Linear Processes with Long Series: Application to Fractional Differencing''. {\slshape Biometrika}, \textbf{83}(3), 603-614.
\bibitem{ni} Nielsen, F.S. (2011). ``Local Whittle Estimation of Multi-Variate Fractionally Integrated Processes''. {\slshape Journal of Time Series Analysis}, \textbf{32}(3), 317-335.
\bibitem{p2} Priestley, M.B. (1981). {\slshape Spectral Analysis and Time Series}. London: Academic Press.
\bibitem{pl} Pumi, G. and Lopes, S.R.C (2013). ``A Semiparametric Estimator for Long-Range Dependent Multivariate Processes''. Under review.
\bibitem{rb1} Robinson, P.M. (1995a). ``Log-Periodogram Regression of Time Series with Long Range Dependence''. {\slshape Annals of Statistics}, \textbf{23}(3), 1048-1072.
\bibitem{rb2} Robinson, P.M. (1995b). ``Gaussian Semiparametric Estimation of Long Range Dependence''. {\slshape Annals of Statistics}, \textbf{23}(5), 1630-1661.
\bibitem{s}  Shimotsu, K. (2007). ``Gaussian Semiparametric Estimation of Multivariate Fractionally Integrated Processes''. {\slshape Journal of Econometrics}, \textbf{137}, 277-310.
\bibitem{sp} Shimotsu, K. and Phillips, P.C.B. (2005). ``Exact Local Whittle Estimation of Fractional Integration''. {\slshape Annals of Statistics}, \textbf{33}(4), 1890-1933.
\bibitem{so}  Sowell, F. (1989). ``Maximum Likelihood Estimation of Fractionally Integrated Time Series Models''. Working Paper, Carnegie-Mellon University.
\bibitem{ts} Tsay, W-J. (2010). ``Maximum Likelihood Estimation of Stationary Multivariate ARFIMA Processes''. {\slshape Journal of Statistical Computation and Simulation}, \textbf{80}(7-8), 729-745.
\bibitem{vel} Velasco, C. (1999). ``Gaussian Semiparametric Estimation of Non-Stationary Time Series''. {\slshape Journal of Time Series Analysis}, \textbf{20}, 87-127.
\bibitem{zyg} Zygmund, A. (2002). {\slshape Trigonometric Series. Vol. I, II}. Cambridge: Cambridge University Press.
\end{thebibliography}
\end{document}